\documentclass[reqno]{amsart}

\usepackage{amsthm}

\usepackage{
amssymb,amsmath,enumerate,stackrel,tabu,mathrsfs,enumitem,verbatim,stmaryrd,
leftidx,bbm, booktabs}
\usepackage{xcolor} 
\usepackage{xr-hyper}
\usepackage{nicefrac,stackrel}
\usepackage[all,cmtip]{xy}
\usepackage[utf8]{inputenc} 
\chardef\cprime"7E 
\usepackage{cases}

\usepackage[T1]{fontenc} 

\usepackage[hyperindex,linktocpage=false]{hyperref}
\hypersetup{
    colorlinks,
    linkcolor={red!50!black},
    citecolor={red!50!black}, 
    urlcolor={red!50!black}, 
    filecolor={red!50!black} 
}

\usepackage{tikz}
\usetikzlibrary{matrix,arrows,cd}

\usepackage[top=1in, left=1in, right=1in, bottom=1in, includeheadfoot]{geometry}

\usepackage{imakeidx}
\makeindex

\usepackage{theoremref}
\definecolor{labelkey}{rgb}{1,0,0}

\makeatletter
\@addtoreset{equation}{section}
\makeatother

\numberwithin{equation}{section}

\theoremstyle{definition}
\newtheorem{Defi}{Definition}[section] \newcommand{\defi}{\begin{Defi}} \newcommand{\xdefi}{\end{Defi}} 
\newtheorem{DefiLemm}[Defi]{Definition and Lemma} \newcommand{\defilemm}{\begin{DefiLemm}} \newcommand{\xdefilemm}{\end{DefiLemm}} 
\newtheorem{Bsp}[Defi]{Example} \newcommand{\exam}{\begin{Bsp}} \newcommand{\xexam}{\end{Bsp}} 
\newtheorem{Syno}[Defi]{Synopsis} \newcommand{\syno}{\begin{Syno}} \newcommand{\xsyno}{\end{Syno}} 
\newtheorem{Bem}[Defi]{Remark} \newcommand{\rema}{\begin{Bem}} \newcommand{\xrema}{\end{Bem}} 
\newtheorem{Notation}[Defi]{Notation} \newcommand{\nota}{\begin{Notation}} \newcommand{\xnota}{\end{Notation}} 

\theoremstyle{plain}
\newtheorem{Theo}[Defi]{Theorem} \newcommand{\theo}{\begin{Theo}} \newcommand{\xtheo}{\end{Theo}} 
\newtheorem{Satz}[Defi]{Proposition} \newcommand{\prop}{\begin{Satz}} \newcommand{\xprop}{\end{Satz}} 
\newtheorem{Lemm}[Defi]{Lemma} \newcommand{\lemm}{\begin{Lemm}} \newcommand{\xlemm}{\end{Lemm}} 
\newtheorem{Coro}[Defi]{Corollary} \newcommand{\coro}{\begin{Coro}} \newcommand{\xcoro}{\end{Coro}}
\newtheorem{Ques}[Defi]{Question} \newcommand{\ques}{\begin{Ques}} \newcommand{\xques}{\end{Ques}}
\newtheorem{Conj}[Defi]{Conjecture} \newcommand{\conj}{\begin{Conj}} \newcommand{\xconj}{\end{Conj}}

\newcommand{\refsect}[1]{Section \ref{sect--#1}}

\newcommand{\eqn}{\begin{equation}} \newcommand{\xeqn}{\end{equation}}
\newcommand{\eqnarr}{\begin{eqnarray*}} \newcommand{\xeqnarr}{\end{eqnarray*}}
\newcommand{\eqnarra}{\begin{eqnarray}} \newcommand{\xeqnarra}{\end{eqnarray}}

\newcommand{\pf}{\begin{proof}} \newcommand{\xpf}{\end{proof}}

\numberwithin{equation}{section}

\newcommand{\nc}{\newcommand}
\nc{\StP}[1]{\cite[\href{http://stacks.math.columbia.edu/tag/#1}{Tag #1}]{StacksProject}}\nc{\StPd}[3]{\cite[Tags~\href{http://stacks.math.columbia.edu/tag/#1}{#1},
\href{http://stacks.math.columbia.edu/tag/#2}{#2},
\href{http://stacks.math.columbia.edu/tag/#3}{#3}]{StacksProject}} 

\nc{\on}{\operatorname}
\nc{\aff}{{\on{aff}}}
\nc{\modi}{{\on{mod}}} 
\nc{\even}{{\on{even}}}
\nc{\odd}{{\on{odd}}}
\nc{\naive}{{\on{naive}}}
\nc{\hofib}{\on{hofib}}
\nc{\Bun}{\on{Bun}}
\nc{\ad}{{\on{ad}}}
\nc{\lft}{{\on{lft}}}

\nc{\Weil}{{\on{Weil}}} 
\nc{\FWeil}{{\on{FWeil}}} 

\nc{\str}{\on{-}}
\nc{\perf}{{\on{perf}}}
\nc{\Rel}{{\on{Pos}}}
\nc{\lan}{\langle}
\nc{\ran}{\rangle}
\nc{\tw}[1]{\langle #1 \rangle} 

\nc{\bbA}{{\mathbb A}} 
\nc{\bbB}{{\mathbb B}}
\nc{\bbC}{{\mathbb C}}
\nc{\bbD}{{\mathbb D}}
\nc{\bbE}{{\mathbb E}}
\nc{\bbF}{{\mathbb F}}
\nc{\bbG}{{\mathbb G}}
\nc{\bbH}{{\mathbb H}}
\nc{\bbI}{{\mathbb I}}
\nc{\bbJ}{{\mathbb J}}
\nc{\bbK}{{\mathbb K}}
\nc{\bbL}{{\mathbb L}}
\nc{\bbM}{{\mathbb M}}
\nc{\bbN}{{\N}} 
\nc{\bbO}{{\mathbb O}}
\nc{\bbP}{{\mathbb P}} 
\nc{\bbQ}{{\mathbb Q}} 
\nc{\bbR}{{\mathbb R}}
\nc{\bbS}{{\mathbb S}}
\nc{\bbT}{{\mathbb T}}
\nc{\bbU}{{\mathbb U}}
\nc{\bbV}{{\mathbb V}}
\nc{\bbW}{{\mathbb W}}
\nc{\bbX}{{\mathbb X}}
\nc{\bbY}{{\mathbb Y}}
\nc{\bbZ}{{\mathbb Z}}

\nc{\calA}{{\mathcal A}}
\nc{\calB}{{\mathcal B}}
\nc{\calC}{{\mathcal C}}
\nc{\calD}{{\mathcal D}}
\nc{\calE}{{\mathcal E}}
\nc{\calF}{{\mathcal F}}
\nc{\calG}{{\mathcal G}}
\nc{\calH}{{\mathcal H}}
\nc{\calI}{{\mathcal I}}
\nc{\calJ}{{\mathcal J}}
\nc{\calK}{{\mathcal K}}
\nc{\calL}{{\mathcal L}}
\nc{\calM}{{\mathcal M}}
\nc{\calN}{{\mathcal N}}
\nc{\calO}{{\mathcal O}}
\nc{\calP}{{\mathcal P}}
\nc{\calQ}{{\mathcal Q}}
\nc{\calR}{{\mathcal R}}
\nc{\calS}{{\mathcal S}}
\nc{\calT}{{\mathcal T}}
\nc{\calU}{{\mathcal U}}
\nc{\calV}{{\mathcal V}}
\nc{\calW}{{\mathcal W}}
\nc{\calX}{{\mathcal X}}
\nc{\calY}{{\mathcal Y}}
\nc{\calZ}{{\mathcal Z}}

\nc{\bone}{{\mathbbm{1}}}

\nc{\Sht}{{\on{Sht}}}
\nc{\Frob}{{\on{Frob}}}
\nc{\Hecke}{{\on{Hecke}}}
\nc{\inv}{{\on{inv}}}
\nc{\Conv}{{\on{Conv}}}
\nc{\triv}{{\on{triv}}}
\nc{\Isom}{{\on{Isom}}}

\nc{\scrB}{{\mathscr{B}}}
\nc{\scrA}{{\mathscr{A}}}
\nc{\bbf}{{\mathbf{f}}}
\nc{\bba}{{\mathbf{a}}}
\nc{\rig}{{\mathrm rig}}

\nc{\al}{\alpha}
\nc{\be}{\beta}
\nc{\ga}{\gamma}
\nc{\la}{\lambda}
\nc{\qcqs}{{\on{qcqs}}}
\nc{\Bmu}{{\mbox{$\raisebox{-0.59ex}{$l$}\hspace{-0.18em}\mu\hspace{-0.88em}\raisebox{-0.98ex}{\scalebox{2}{$\color{white}.$}}\hspace{-0.416em}\raisebox{+0.88ex}{$\color{white}.$}\hspace{0.46em}$}{}}}

\nc{\pot}[1]{ [\hspace{-0,5mm}[ {#1} ]\hspace{-0,5mm}] }
\nc{\rpot}[1]{ (\hspace{-0,7mm}( {#1} )\hspace{-0,7mm}) }

\nc{\defined}{\hspace{0.1cm}\stackrel{\text{\tiny \rm def}}{=}\hspace{0.1cm}}
\nc{\co}{\colon}
\nc{\specto}{{\leadsto}}

\newcommand{\res}{\operatorname{res}} 

\def\Gm{\mathbf {G}_\mathrm m} 
\def\SL{\mathrm {SL}} 
\def\PGL{\mathrm {PGL}} 

\font\tencyr=wncyr10
\font\sevencyr=wncyr7
\font\fivecyr=wncyr5
\def\To#1#2{\mathop{\count0=#1 \loop\ifnum\count0>0 \smash-\mkern-7mu \advance\count0 -1 \repeat \mathord\rightarrow}\limits^{#2}} 

\def\Gr{\mathop{\rm Gr}\nolimits} 
\def\Fl{\mathop{\rm Fl}\nolimits} 
\def\Sht{\mathop{\rm Sht}\nolimits} 

\definecolor{hellgrau}{RGB}{200,200,200} 
\definecolor{dunkelgrau}{RGB}{160,160,160} 
\definecolor{hellblau}{RGB}{194, 215, 249} %
\definecolor{dunkelblau}{RGB}{68, 128, 226} %

\def\Z{{\mathbb{Z}}} 
\def\N{{\bf N}} 
\def\Q{{\mathbb Q}} 
\def\Gm{\mathbf {G}_\mathrm m} 

\def\Gal{{\rm Gal}} 
\def\der{{\rm der}} 

\def\Spec{\mathop{\rm Spec}} 

\def\sbuildrel#1\over#2{\mathrel{\smash{\mathop{\kern0pt #2}\limits^{#1}}}}

\let\Gets\longleftarrow
\let\x\times

\newcommand{\lr}{\longrightarrow}

\newcommand{\SO}{\mathrm{SO}}
\newcommand{\PSO}{\mathrm{PSO}}
\newcommand{\PU}{\rm{PU}}
\nc{\supp}{\mathrm{supp}}
\nc{\cw}[1]{\omega_{#1}^\vee}
\nc{\fraka}{{\mathfrak{a}}}
\nc{\frakh}{{\mathfrak{h}}}
\nc{\frakz}{{\mathfrak{z}}}

\nc{\Sv}{\mathcal{S}} 
\nc{\BTB}{\mathscr{B}} 
\nc{\BTA}{\mathscr{A}} 

\DeclareMathOperator{\Res}{Res}

\begin{document}

\title[]{Normality of Schubert varieties in \\ affine Grassmannians II: 
	\\ The tamely ramified case}
\author[P.\,Bieker]{Patrick Bieker}

\address{ Fakultät für Mathematik, Universität Bielefeld, Postfach 100 131, 33501 Bielefeld, Germany}
\email{pbieker@math.uni-bielefeld.de}

\thanks{The author acknowledges support by the Deutsche Forschungsgemeinschaft (DFG, German Research Foundation), project number 520675682, as well as by the DFG through the TRR 358 \textit{Integral Structures in Geometry and Representation Theory}, project number 49139240.}

\maketitle

\begin{abstract} 
	We prove a criterion for the normality of Schubert varieties in twisted affine Grassmannians in terms of the order of the algebraic fundamental group of a certain Levi subgroup, in particular in small positive characteristic.
	
	As an application, we obtain a similar normality criterion for local models in both equal and mixed characteristic. 
	In particular, we give a classification of normal Pappas-Zhu local models at absolutely special level as well as for adjoint groups of rank 1.
\end{abstract}


\thispagestyle{empty}

\section{Introduction}

Schubert varieties in (partial) affine flag varieties play an important role in both geometric representation theory and arithmetic geometry, for example via their close relation to local models of Shimura varieties (in mixed characteristic) or moduli spaces of shtukas (in positive characteristic).
One of the fundamental questions is the normality of Schubert varieties for a reductive group $G$ over a field of Laurent series $k \rpot{t}$ and local models for reductive groups $\bbG$ over non-archimedean local fields $\bbF$ (either in mixed or positive characteristic).

When the characteristic of the ground field is zero or big enough in the sense that it does not divide the order of the fundamental group $\pi_1(G^{\rm der})$, Faltings \cite{Faltings:Loops}, Pappas--Rapoport \cite{PappasRapoport:Twisted} and Lourenço \cite{Lourenco:Normality} prove the normality of Schubert varieties. Similarly, Zhu \cite{Zhu:Coherence}, Pappas-Zhu \cite{PappasZhu:Kottwitz}, Levin \cite{Levin:Weil}, Richarz \cite{Richarz:AffGrass}, Haines-Richarz \cite{HainesRicharz:Normality} and Lourenco \cite{Lourenco:Grassmannian} construct and prove normality of local models whenever the residue characteristic of $\bbF$ does not divide the order of $\pi_1(\bbG^{\der})$. 

In small positive characteristic, when ${\rm char}(k)\mid \#\pi_1(G^{\rm der})$, Haines--Louren\c{c}o--Richarz \cite{HainesLourencoRicharz:Normality} produce examples of non-normal Schubert varieties and show that in the absolutley almost simple case  almost all Schubert varieties in the affine flag variety are non-normal. As a consequence also the local models are non-normal in this case.
Subsequently, the author together with Richarz \cite{BiekerRicharz:normality} classify all (non-)normal Schubert varieties in affine Grassmannians for split semisimple and almost simple groups.

The main goal of the present manuscript is to give a criterion for the normality of Schubert varieties in the affine Grassmannian $\Gr_G$ for many reductive groups $G$, including all split groups as well as all adjoint groups in characteristic ${\rm char}(k) \geq 5$ or that split over a tamely ramified extension in charactersitic 2 or 3, in terms of the order of the fundamental group of a certain Levi subgroup.
We use the criterion to extend the classification result of \cite{BiekerRicharz:normality} to tamely ramified and not necessarily semisimple groups.
As an application we give an analogous criterion for the normality of local models and classify the normal local models at absolutely special level as well as for groups of semisimple rank 1. 

\subsection{Results} 
Let $k$ be an algebraically closed field, and let $F \defined k \rpot{t} $ be the field of formal Laurent series over $k$ with absolute Galois group $I = \Gal(F^{\rm sep}/F)$.
Let $G$ be a reductive group over $F$, then $G$ is automatically quasi-split.
To formulate our main theorem, we impose the condition on $G$ that 
\eqn
G \simeq \prod_{j \in J} {\rm Res}_{F_j/F} \widetilde{G}_j
\label{eqn:cond-G}
\xeqn
for a finite set $J$, finite separable extensions $F_j/F$ and reductive groups $\widetilde{G}_j$ over $F_j$ that split over a tamely ramified extension $F_j'$ of $F_j$.
In particular, we allow for example all tamely ramified groups and also all adjoint groups in characteristic ${\rm char} (k) \geq 5$.

\subsubsection{A normality criterion}
For a special vertex $x$ in the Bruhat-Tits building $\BTB(G,F)$ there is the associated affine Grassmannian $\Gr_{G,x}$.
Fix a maximally split maximal torus $T$ contained in a Borel subgroup $B$ in $G$ such that $x$ is contained in the apartment corresponding to the maximal split torus $S$ in $T$.
Associated to each $B$-dominant cocharacter $\bar \mu\in X_*(T)^+_I$ there is the Schubert variety $\Gr_{G,x, \leq \bar \mu}$, defined as a certain orbit closure in the affine Grassmannian $\Gr_{G, x}$, see \refsect{Schubert-schemes}. 

On the set of dominant cocharacters $X_*(T)^+_I$ there is the Bruhat order which is explicitly given by $\bar \la \leq \bar \mu$ if and only if $\bar \mu - \bar \la$ is a linear combination with non-negative coefficients of (\'echellonage) coroots.
The set of coroots that appear with positive coefficient is called the support of $\bar \mu - \bar \la$. 
For every $\bar \mu \in X_*(T)^+_I$ exists a unique minuscule $\bar \la \in X_*(T)^+_I$ with $\bar \la \leq \bar \mu$. We denote by $M_{\bar \mu} \subseteq G$ the Levi subgroup corresponding to the subset ${\rm supp}(\bar \mu - \bar \la)$ of positive relative coroots.

\theo[{cf.\@ \thref{thm:main}}]
\thlabel{thm:normality-criterion}
Let $G$ be a reductive group that satisfies \eqref{eqn:cond-G} and let $x \in \BTB(G,F)$ be a special vertex. 
Let  $\bar \mu \in X_*(T)^+_I$ be a dominant cocharacter. If
\eqn 
{\rm char}(k)\nmid \#\pi_1(M_{\bar \mu}^{\der}),
\label{eq:char-pi1}
\xeqn
then the Schubert variety $\Gr_{G,x, \leq \bar \mu}$ is normal.
Moreover, if $x$ is an absolutely special vertex and $G$ is split or adjoint, then \eqref{eq:char-pi1} is also necessary for the normality of $\Gr_{G,x, \leq \bar \mu}$.
\xtheo

The conclusion in the last sentence is false for arbitrary (even almost simple) tamely ramified groups as the quasi-minuscule Schubert variety for $G = ({\rm Res}_{F'/F} \SL_2) / \Bmu_2$ is normal for $F$ of characteristic 2 and $F'/F$ an extension of odd degree, compare \thref{prop:res-sl2-normal}.
Moreover, the fundamental group appearing in the condition \eqref{eq:char-pi1} can be computed effectively from the simple coroots of $M_{\bar \mu}$ as 
\eqn
\pi_1(M^{\der}_{\bar \mu}) = (X_*(T) \cap \bbQ\Phi_{\rm abs}^\vee(M_{\bar \mu},T))/\bbZ\Phi_{\rm abs}^\vee(M_{\bar \mu},T)
\label{eqn:order-fundamental-group}
\xeqn
by \cite[Lemma 5.3]{BiekerRicharz:normality}. In many cases it even suffices to determine the connection index of the type of $M_{\bar \mu}$ to deduce that \eqref{eq:char-pi1} is satisfied.

As an application, we deduce a similar normality criterion for local models. 
Let now $\bbG$ be a connected reductive group over a non-archimedean local field $\bbF$ (either in characteristic 0 or positive characteristic) with residue field $\kappa$. 
We denote by $\breve{\bbF}$ the completion of the maximal unramified extension of $\bbF$. 
Assume that $\bbG$ satisfies the analogue of condition \eqref{eqn:cond-G}, i.e. $\bbG$ is a product of restrictions of scalars of reductive groups that split over a tamely ramified extension of $\bbF$. 
For a conjugacy class $\{\mu\}$ of (geometric) cocharacters of $\bbG$ and a facet $\bbf$ in the Bruhat-Tits building $\BTB(\bbG, \bbF)$ we have the local model
$$ \calM^{\rm loc} = \calM^{\rm loc}(\bbG, \{\mu\}, \calG_{\bbf})$$
defined over the ring of integers of the reflex field
as constructed by Pappas-Zhu \cite{PappasZhu:Kottwitz} and Levin \cite{Levin:Weil} in mixed characteristic as well as Zhu \cite{Zhu:Coherence} and Richarz \cite{Richarz:AffGrass} in positive characteristic. 
As a corollary to \thref{thm:normality-criterion} together with a result of Haines-Richarz \cite{HainesRicharz:Normality} we obtain a similar normality criterion for local models.
\theo[{cf.\@ \thref{theo:normality-loc-mod}}]
\thlabel{thm:normality-criterion-locmod}
	Let $\bbG$ be a reductive group over $\bbF$ that satisfies Condition \eqref{eqn:cond-G} and let $\mathbf f \in \BTB(\bbG, \bbF)$ be a facet that contains a special vertex in its closure.
	Then the local model $\calM^{\rm loc}$ is normal if
	\eqn 
	{\rm char}(\kappa)\nmid \#\pi_1(M_{\bar \mu}^{\der}).
	\label{eqn:char-pi1-locMod}
	\xeqn
	Moreover, if $\mathbf f = \mathbf 0$ is an absolutely special vertex and $\bbG$ is adjoint or splits over an unramified extension of $\bbF$, then \eqref{eqn:char-pi1-locMod} is also necessary for the normality of $\calM^{\rm loc}$.	
\xtheo

\subsubsection{Classification of normal Schubert varieties}
We use our normality criteria to give a full classification of normal Schubert varieties (respectively normal local models) for almost simple groups at absolutely special level generalizing \cite[Theorem 1.1]{BiekerRicharz:normality} to tamely ramified and not necessarily semisimple groups. 
Following \cite{HainesLourencoRicharz:Normality} we use the convention that a reductive group $G$ is called \emph{absolutely almost simple} if its absolute Dynkin diagram is connected. The group $G$ is not assumed to be semisimple.

Let us fix an absolutely special vertex $\mathbf 0 \in B(G,F)$. 
Let $P^{\vee}$ be the coweight lattice of the \'echellonage roots system for $G$, which is a reduced root system in $(X_*(T)^+_I)_{\bbR}$, compare \refsect{echellonage-roots}, and let $n\in\Z_{\geq 0}$ be its rank. 
We denote by $\cw1,\ldots,\cw n$ the $\Z$-basis of $P^\vee$ of fundamental coweights (corresponding to the choice of $T \subseteq B$ from above) in the notation of \cite[Tables]{Bourbaki:Lie456}.
We have a map of partially ordered monoids 
\eqn
X_*(T)_I^+ \to (P_G^\vee)^+
\label{eqn:map-adjoint}
\xeqn
and $\cw1,\ldots,\cw n$ form a monoid basis of $(P_G^\vee)^+$.
We note that the isomorphism type of $\Gr_{G,\mathbf 0, \leq \bar \mu}$ only depends on the image $\bar \mu^{\rm ad}$ of $\bar \mu$ under \eqref{eqn:map-adjoint}.
Together with the normality result when ${\rm char}(k)\nmid \#\pi_1(G^{\rm der})$ the following theorem classifies all normal Schubert varieties in $\Gr_{G, \mathbf 0}$:

\theo[{\refsect{classification}, cf. also \cite[Theorem 1.1]{BiekerRicharz:normality}}]
\thlabel{classification-normal-schubert-varieties}
Let $G$ be an absolutely almost simple (but not necessarily semisimple) reductive group that splits over a tamely ramified extension of $F$. 
Let $\mathbf 0 \in \BTB(G,F)$ be an absolutely special vertex.
Let $\bar \mu\in X_*(T)^+_I$ be dominant cocharacter. Assume that ${\rm char}(k)\mid \#\pi_1(G^{\rm der})$.
Then, $\Gr_{G,\mathbf 0, \leq \bar \mu}$ is normal if and only if the pair $(G, \bar \mu)$ belongs to the following list:
\begin{itemize}
	\item any $G$, and $ \bar \mu$ minuscule; 
	\item $G$ of type $A_n$ for $n\geq 2$ \textup{(}so ${\rm char}(k)\mid n+1$\textup{)}, and $\bar \mu^{\rm ad}\leq d\cw 1$ or $\bar \mu^{\rm ad}\leq d\cw n$ for some $d\in\{2,\ldots,n\}$;
	\item $G$ of type $D_n$ for $n \geq 4$ \textup{(}so ${\rm char}(k) = 2$\textup{)} and 
	\begin{itemize}
		\item 
		$\bar \mu^{\rm ad}\in\{\cw1+\cw{n-1}, \cw1+\cw{n}\}$ and $G^{\rm der} \simeq \PSO_{2n,k}$ for $n$ odd or $G^{\der} \simeq \SO_{2n, k}$ (for all $n \geq 4$); 
		\item 
		$G^{\rm der} \not \simeq \SO_{2n,k}, \PSO_{2n,k}$, in particular $n \geq 6$ even, and $\bar \mu^{\rm ad} \leq \cw{n-1} + \cw n$; 
		\item $G^{\rm der} \not \simeq \SO_{2n,k}, \PSO_{2n,k}$ for $n \equiv 2 \mod 4$, and $\bar \mu^{\rm ad} = \cw 1 + \cw{n-1}$;
		\item $G^{\rm der} \not \simeq \SO_{2n,k}, \PSO_{2n,k}$ for $4 \mid n$, and $\bar \mu^{\rm ad} = \cw 1 + \cw{n}$;
	\end{itemize}
	\item $G$ of type $E_6$ \textup{(}so ${\rm char}(k)=3$\textup{)}, and $\bar \mu^{\rm ad}\in\{2 \cw 1, \cw  3, \cw 5, 2\cw 6\}$;
	\item $G$ of type $E_7$ \textup{(}so ${\rm char}(k)=2$\textup{)}, and $\bar \mu^{\rm ad}=\cw 2$.
	\item $G$ of type $B\str C_n$ \textup{(}so ${\rm char}(k)\mid 2n$ and  ${\rm char}(k) \neq 2$\textup{)}, and $\bar \mu^{\rm ad} = \cw{m}$ for some odd $1 \leq m \leq n$.
\end{itemize}
\xtheo

The (non-)normality of Schubert varieties for absolutely almost simple groups outside type $D_n$ does not depend on the isomorphism type of $G$ nor on the characteristic of $k$ as long as ${\rm char}(k) \mid \#\pi_1(G^{\rm der})$.
Recall that in type $D_n$ for $n \geq 6$ and $n$ even apart from ${\rm Spin}_{2n,k}, \SO_{2n,k}$ and $\PSO_{2n,k}$ there is a fourth isomorphism class of semisimple groups, this is the group $G \not \simeq \SO_{2n,k}, \PSO_{2n,k}$ appearing in the statement of the theorem. 
We identify its root datum with the data in \cite{Bourbaki:Lie456} in such a way that $\cw n$ is a cocharacter of this group.

In a similar fashion we classify pairs $(\bbG, \{\mu\})$, where $\bbG$ is defined over a non-archimedean local field $\bbF$ as above, such that the local model at absolutely special level is normal.
To $\bbG$ we associate an equal-characteristic analogue $G$ of the same type defined over $F \defined k \rpot{t}$ where $k = \kappa^{\rm alg}$.
In particular, when $\rm{char}(\bbF)>0$, then $G = \bbG \otimes_{\bbF} F$.

\theo[{\thref{theo:normality-loc-mod}}]
Let $\bbG$ be an absolutely almost simple group that splits over a tamely ramified extension of $\bbF$ and fix an absolutely special vertex $\mathbf 0$. Let $\{\mu\}$ be a conjugacy class of geometric cocharacters of $\bbG$ such that $\{\mu\}$ is not central. 
Moreover, assume ${\rm char}(\kappa) \mid \# \pi_1(\bbG^{\der})$.
Then the local model $\calM^{\rm loc} = \calM^{\rm loc}(\bbG, \{\mu\}, \calG_{\mathbf 0})$ is normal if and only if $(\bbG, \{\mu\})$ belongs to the following list:
\begin{enumerate}
	\item $\bbG$ splits over $\breve \bbF$ and $(G, \{\mu\})$ appears in the list in \thref{classification-normal-schubert-varieties}.
	\item $\bbG_{\breve \bbF}$ is of type $B\str C_n$ \textup{(}so ${\rm char}(\kappa)\mid 2n$ and  ${\rm char}(\kappa) \neq 2$\textup{)}, and $\mu = \cw{2i-1}$ for some $1 \leq i \leq n$. 
\end{enumerate}
\xtheo

We also finish the classification of Schubert varieties in the affine flag variety for adjoint groups of semisimple $F$-rank 1, compare \cite[Corollary 6.8]{HainesLourencoRicharz:Normality} for the case $G \simeq \PGL_2$ in characteristic 2, and \thref{prop:classification-PU3} for the case $G \simeq \PU_3$ in characteristic 3.
As a consequence, we can classify the normal local models for arbitrary parahoric level in this case.
\theo[{\thref{thm:class-locmod-ss1}}]
Let $(\bbG, \{\mu\}, \calG_{\bbf})$ be an LM-triple such that $\bbG_{\breve \bbF}$ is adjoint of $\breve \bbF$-rank 1 and its absolutely simple factor splits over a tamely ramified extension.
Assume moreover ${\rm char}(\kappa) \mid \pi_1(\bbG)$. 
Then $\calM^{\rm loc}(\bbG, \{\mu\}, \calG_\bbf)$ is normal if and only if the LM-triple appears in the following list.
\begin{enumerate}
	\item $\bbG_{\breve \bbF} \simeq {\rm Res}_{\breve \bbF'/\breve \bbF} \PGL_{2, \breve{\bbF}'}$ (so ${\rm char}(\kappa) = 2$) and 
	\begin{enumerate}
		\item $\mu = (0, \ldots, 0, \widetilde{\mu}_j, 0, \ldots,0)$ with $\widetilde{\mu}_j \in X_*(\widetilde{T})^+$ minuscule (for arbitrary $\calG_\bbf$ and $\bbF$), or
		\item $\mu = (0, \ldots, 0, \widetilde{\mu}_j, 0, \ldots,0, \widetilde{\mu}_{j'}, 0, \ldots,0)$ with $\widetilde{\mu}_j, \widetilde{\mu}_{j'} \in X_*(\widetilde{T})^+$ minuscule and $\calG_\bbf= \calG_\bba$ is an Iwahori group scheme, or
		\item $\mu = (0, \ldots, 0, \widetilde{\mu}_j, 0, \ldots,0)$ with $\widetilde{\mu}_j \in X_*(\widetilde{T})^+$ quasi-minuscule, ${\rm char}(\bbF) = 0$ and $\calG_\bbf = \calG_\bba$ is an Iwahori group scheme.
	\end{enumerate}
	\item $\bbG_{\breve \bbF}  \simeq {\rm Res}_{\breve \bbF'/\breve \bbF} \PU_{3, \breve{\bbF}'}$ (so ${\rm char}(\kappa) = 3$) and
	\begin{enumerate}
		\item $\mu = 0$, or
		\item $\mu = (0, \ldots, 0, \widetilde{\mu}_j, 0, \ldots,0)$ with $\widetilde \mu_j \in \{\cw 1, \cw 2\}$ is minuscule and 
		\begin{itemize}
			\item $\bbf = x$ is a special but not absolutely special vertex or 
			\item $\bbf = \bba$ is an alcove.
		\end{itemize}
	\end{enumerate}
\end{enumerate}
\xtheo
In the formulation of the theorem we use that the geometric cocharacters of a maximal torus ${T} = {\rm Res}_{F'/F} \widetilde{T}$ of $G = {\rm Res}_{F'/F} \widetilde G$, where $\widetilde{T}$ is a maximal torus of $\widetilde{G}$, decompose as 
$$X_*(T) \cong \bigoplus_{[F'\colon F]} X_*(\widetilde{T}).$$ 

The (non-)normality of Schubert varieties and local models can then via the local model theorems translated into corresponding statements concerning the geometry of integral models of Shimura varieties and moduli spaces of shtukas.
For example, \thref{thm:normality-criterion} (even only its split case) gives a similar criterion for the normality of generic fibers of moduli spaces of shtukas. 
Similarly, the criterion and classification for the normality of local models implies corresponding results for the full moduli space of moduli spaces of shtukas with parahoric level structure in small positive characteristic.

\subsection{Proof strategy}
For the proof of our main theorems we use and generalize the methods of \cite{BiekerRicharz:normality} and \cite{HainesLourencoRicharz:Normality}. 
In \refsect{transversal-slice}, we prove a twisted version of the Levi lemma of \cite{MalkinOstrikVybornov:MinimalDegenerations} in our setting in positive characteristic, generalizing also similar results in \cite{BiekerRicharz:normality}, \cite{Richarz:Diplom} and \cite{BessonHong:SmoothLocus}.  
The Levi lemma shows in particular, that the Schubert variety $\Gr_{G,x, \leq \bar \mu}$ is smoothly equivalent to the Schubert variety $\Gr_{M_{\bar \mu}, x, \leq \bar \mu}$, where $M_{\bar \mu} \subseteq G$ is the Levi subgroup defined above. 
Together with the normality of Schubert varieties when ${\rm char}(k)$ does not divide the order of the fundamental group, c.f. \cite{PappasRapoport:Twisted}, this already shows one implication of \thref{thm:normality-criterion}.

In order to show the non-normality we generalize the tangent space calculations in the proof of the non-normality of the quasi-minuscule Schubert variety in \cite[Corollary 6.2]{HainesLourencoRicharz:Normality}.
In \refsect{quasi-minuscule} we show that for not necessarily absolutely almost simple groups the non-normality of the Schubert variety corresponding to the cocharacter that is quasi-minuscule in every absolutely almost simple factor.
We then in Section \ref{sect--normality-criterion} deduce \thref{thm:normality-criterion} by combining the results of the previous two sections.

In order to deduce the classification result \thref{classification-normal-schubert-varieties} from our normality criterion we determine in \refsect{classification} all dominant cocharacters that are not directly covered by the non-normality result  \thref{prop:intro}. 
These are only finitely many (up to the center), and in each case we compute the order of the fundamental group of the corresponding Levi subgroup.
This then covers all the cases in the classification that are not already covered in the split semisimple case in \cite[Theorem 1.1]{BiekerRicharz:normality}.

In \refsect{general-level} we deduce partial results for more general (i.e., not necessarily absolutely special) level structure. 
In particular, we give the classification of normal Schubert varieties in the (full) affine flag variety for $\PU_3$ and classify all normal Iwahori-Schubert varieties in the affine Grassmannian for $\PGL_n$. 
In particular, we produce many new examples of normal Schubert varieties in the affine flag variety for $\PGL_n$. 

In the \refsect{loc-model} we apply our results to the geometry of local models. 
We make use of results of \cite{HainesRicharz:Normality} that the local model is normal if and only if all Schubert varieties that appear in the generic fiber or the admissible locus in the special fiber are normal. 
But using the results from before, we can answer these questions in many cases.

\medskip
\noindent\textbf{Acknowledgements.} 
It is a pleasure to thank R{\i}zacan \c{C}ilo\u{g}lu, Manuel Hoff and Timo Richarz for many helpful conversations surrounding this work.


\section{Root systems and Weyl groups}\label{sect--Schubert-schemes}
We recall some relevant definitions and facts related to Weyl groups and root systems. 
Let $k$ be an algebraically closed field and let $F = k \rpot{t}$ be the field of Laurent series over $k$ with ring of integers $\calO = k \pot{t}$. We fix a separable closure $F^{\rm sep}$ of $F$ and denote by $I = \Gal(F^{\rm sep}/F)$ its absolute Galois group.
Let $G$ be a connected reductive group over $F$, then $G$ is automatically quasi-split.
We denote by $F'$ the minimal splitting field of $G$, it is a totally ramified cyclic Galois extension of $F$ of degree $e$. 

Let $\bbf$ be a facet in the (extended) Bruhat-Tits building $\BTB(G,F)$ and denote by $\calG_{\bbf}$ the corresponding parahoric group scheme, it is a smooth affine group scheme with connected fibres over $\calO$ and generic fiber $G$.
Fix a maximal $F$-split torus $S$ in $G$ such that $\bbf$ is contained in the apartment $\BTA(G,S,F) \subseteq \BTB(G,F)$ corresponding to $S$.
We denote by $T$ the centralizer of $S$ in $G$, it is a maximal $F$-torus in $G$ as $G$ is quasi-split. We fix a Borel subgroup $B \subseteq G$ containing $T$. 

Moreover, let $G^{\rm{der}} \subseteq G$ be its derived subgroup with simply connected cover $G^{\rm{sc}} \to G^{\rm{der}}$ let $G^{\rm{ad}}$ be its adjoint quotient.
Let 
\[
	T^{\rm{der}}, T^{\rm{sc}}, T^{\rm{ad}}, S^{\rm{der}}, S^{\rm{sc}} \text{ and } S^{\rm{ad}}
\]
be the (pre-)images of $T$ and $S$ in $G^{\rm{der}}, G^{\rm{sc}}$ and $G^{\rm{ad}}$, respectively. Then $T^*$ and $S^*$ are maximal ($F$-split) tori in $G^*$ for $* \in \{\rm{der}, \rm{sc}, \rm{ad}\}$.


\subsection{Relative and absolute root datum}
We denote by $X^*(T)$ (respectively $X_*(T)$) the \emph{geometric characters} (respectively \emph{geometric cocharacters}) of $T$.
We denote by $\lan\str,\str\ran$ the natural pairing between the characters $X^*(T)$ and cocharacters $X_*(T)$.
Let $\Phi_{\rm{abs}} = \Phi(G,T) \subset X^*(T)$ and $\Phi^\vee_{\rm{abs}} = \Phi^\vee(G,T) \subset X_*(T)$ be the set of \emph{absolute (co-)roots} of $(G,T)$.
In particular, $\Z \Phi_{\rm{abs}} = X^*(T^{\rm{sc}})$.

Let $\Phi_{\rm{abs}}^+ \subset \Phi_{\rm{abs}}$ be the set of $B$-positive roots and $\Delta_{{\rm abs}} \subseteq \Phi_{\rm{abs}}^+$ a basis of positive roots. 
We denote by $\Delta_{{\rm abs}}^\vee$ the corresponding set of simple coroots.
The sum of the $B$-positive roots is denoted $2\rho\in X^*(T)$.
Let $(X_*(T)^+,\leq)$ be the partially ordered monoid of $B$-dominant cocharacters equipped with the Bruhat partial order, i.e., 
for $\la,\mu \in X_*(T)^+$ one has $\la\leq\mu$ if and only if $\mu-\la$ is a sum of positive coroots with $\Z_{\geq 0}$-coefficients.

The group $I$ acts on the lattices $X^*(T), X_*(T), \Z\Phi_{\rm{abs}}$ and $\Z\Phi^\vee_{\rm{abs}}$ through its finite quotient $\Gal(F'/F)$. 
For any $I$-module $M$ we denote by $M_I$ its coinvariants. 
When $M$ is a free abelian group that admits an $I$-stable basis, also the group of coinvariants $M_I$ is a free abelian group. 
When $G$ is simply connected or adjoint the (co-)character lattices admit $I$-stable bases given by the (co-)roots or fundamental (co-)weights, respectively.
The maps
\[
	T^{\rm{sc}} \to T^{\rm{der}} \to T \to T^{\rm{ad}}
\]
induce $I$-equivariant maps on the corresponding cocharacter lattices
\[
	X_*(T^{\rm{sc}}) \hookrightarrow X_*(T^{\rm{der}}) \hookrightarrow X_*(T) \to X_*(T^{\rm{ad}}),
\]
which descend to maps on coinvariants
\[
	X_*(T^{\rm{sc}})_I \hookrightarrow X_*(T^{\rm{der}})_I \hookrightarrow X_*(T)_I \to X_*(T^{\rm{ad}})_I.
\]
The quotient of the cocharacter lattice by the coroot lattice is called the \emph{algebraic fundamental group} of $G$, denoted by
\[
	\pi_1(G) \defined X_*(T)/(\Z\Phi^\vee_{\rm{abs}}).
\]
It is well-known that $\pi_1(G^{\rm{der}})$ is a finite abelian group.
The action of $I$ descends to an action on $\pi_1(G)$ and induces an isomorphism 
\[
	\pi_1(G)_I \cong X_*(T)_I/(\Z\Phi^\vee_{\rm{abs}})_I
\]
by \cite[Lemma 2.5]{AcharLourencoRicharzRiche:Modular}.
Moreover, there is the so-called \emph{Kottwitz homomorphism}
\[
	\kappa_G \colon G(F) \to \pi_1(G)_I
\]
as defined in \cite[Section 7]{Kottwitz:Isocrystals2}.

Similarly, let $X^*(S)$ and $X_*(S)$ be the lattice of cocharacters of $S$ containing the set of relative (co-)roots $\Phi = \Phi(G,S) \subset X^*(S)$ and $\Phi^\vee(G,S) \subset X_*(S)$, respectively.  
The inclusion $S \subseteq T$ induces a projection $X^*(T) \to X^*(S)$ which factors through $X^*(T)_I$ and identifies $X^*(S)$ with the maximal torsion-free quotient of $X^*(T)_I$. Moreover, the set of absolute roots $\Phi_{\rm{abs}}$ is mapped onto $\Phi$, and the image $\Phi^+$ of $\Phi^+_{\rm{abs}}$ is a set of positive roots for $\Phi$. Let $\Delta \subseteq \Phi^+$ be a basis of ($B$-positive) roots.

The submonoid of \emph{($B$-)dominant} elements of $X_*(T)_I$ is given by 
\eqn
	X_*(T)^+_I \defined \left\{ \bar \mu \in X_*(T)_I \colon \forall \alpha \in \Phi^+, \langle \bar \mu, \alpha \rangle \geq 0 \right\}.
\xeqn 
 
\subsection{Iwahori-Weyl group}
Let $W_0 = W(G,S) = N_G(T)(F)/T(F)$ be the finite Weyl group of $(G,S)$, and let 
\eqn
	W = W(G,S) \defined N_G(T)(F)/(T(F) \cap \ker \kappa_T)  
\xeqn
be the associated \emph{Iwahori-Weyl group}.
Then $W$ acts on the apartment $\BTA(G,S,F)$ by affine transformations.
For tori there is an isomorphism
\begin{align*}
	X_*(T)_I & \xrightarrow{\cong} W(T,S) = T(F)/(T(F) \cap \ker \kappa_T) \\
	\bar \mu & \mapsto \bar \mu(t). 
\end{align*}
Recall that more explicitly when $T$ is split an explicit lift of $\mu(t)$ to $T(F)$ is given by the image of $t \in F^{\times}$ under $\mu \colon \Gm(F) \to T(F)$, justifying the notation.
The inclusion $T \subset G$ induces an exact sequence
\eqn
	0 \to X_*(T)_I \to W \to W_0 \to 1.
	\label{eqn:extension-Weyl-grp}
\xeqn 

For any facet $\bbf$ of $\BTA(G,S,F)$ let $W_{\bbf} \subset W$ be the subgroup generated by the  reflections along root hyperplanes containing $\bbf$.  
Note that $W_{\bbf}$ is a finite Coxeter group. 
Moreover, for a special vertex $x \in \BTA(G,S,F)$ (which always exists as $G$ is quasi-split) the natural map $W_x \hookrightarrow W \twoheadrightarrow W_0$ is an isomorphism, and thus induces a splitting
\eqn
	W \cong X_*(T)_I \rtimes W_0
	\label{eqn:ext-affine-Weyl-grp}
\xeqn
identifying the Iwahori-Weyl group with the \emph{extended affine Weyl group}.

The set of affine hyperplanes in $\BTA(G,S,F)$ through the special point $x$ determines chambers in $\BTA(G,S,F)$ and our choice of Borel subgroup $B$ singles out a positive one. 
As base alcove $\bba$ we fix the alcove in the negative chamber that contains $x$ in its closure.
Let $\Omega_{\bba} \subset W$ be the stabilizer of $\bba$.
Let $W_{\rm{aff}} \subset W$ be the subgroup generated by reflections along the walls of $\bba$, it is called the \emph{affine Weyl group}. Note that $W_{\rm{aff}} = W(G^{\rm{sc}}, S^{\rm{sc}})$ is the Iwahori-Weyl group of the simply connected cover of $G$. 
Moreover, $W$ acts simply transitively on the set of alcoves in $\BTA(G,S,F)$ and so we have
\eqn
	W \cong W_{\rm{aff}} \rtimes \Omega_{\bba}.
	\label{eq:W-Waff}
\xeqn
Additionally, the composition $X_*(T)_I \hookrightarrow W \twoheadrightarrow \Omega_{\bba}$ induces an isomorphism
\eqn
	\pi_1(G)_I \xrightarrow{\cong} \Omega_{\bba}.
	\label{eq:pi1-Omega}
\xeqn
By construction, the affine Weyl group is endowed with the structure of a Coxeter group, yielding a quasi-Coxeter structure on $W$ using the splitting \eqref{eq:W-Waff}. 
In particular, we have a length function $\ell$ and the Bruhat order $\leq$ on $W$.
The elements of length 0 in $W$ are by definition given by $\Omega_{\bba} \subseteq W$.

For a pair of facets $\bbf, \bbf'$ in $\BTA(G,S,F)$ we also consider the induced order on double quotients $W_{\bbf'} \backslash W / W_{\bbf}$. 
Recall that by the Bruhat-decomposition we have  $W_{\bbf'} \backslash W / W_{\bbf} \cong \calG_{\bbf'}(\calO) \backslash G(F)/ \calG_{\bbf}(\calO)$.
For a choice of special vertex we then get using \eqref{eqn:ext-affine-Weyl-grp} that the natural maps  
\[
	X_*(T)_I \xrightarrow{\cong} W/W_x 
	\qquad \text{ and } \qquad 
	X_*(T)^+_I \xrightarrow{\cong} W_x \backslash W  / W_x
\]
are isomorphisms.
In particular, we have the (quotient) Bruhat order $\leq$ on $X_*(T)_I$ and $X_*(T)^+_I$. 

Let us also record the following fact.
\lemm $
		\Omega_{\bba} \cap X_*(T)_I = 
		\ker( X_*(T)_I \to X_*(T^{\rm{ad}})_I ) = \ker(W \to W^{\rm ad})
	$
	\thlabel{lem:central-cochars}
\xlemm
\pf
	For the first equality we note that the translation $\bar{\mu}(t)$ acts only in central directions on $\BTA(G,S,F)$, and in particular preserves the alcove $\bba$ (or equivalently all alcoves) if and only if $\bar \mu$ pairs trivially with all (relative) roots.
	As the map $X_*(T)_I \to X_*(T^{\rm{ad}})_I$ preserves the pairing by construction, this gives the first equality. 
	
	The second equality follows from comparing \eqref{eqn:extension-Weyl-grp} for $G$ and $G^{\rm ad}$. 
\xpf

\subsection{Product decompositions}
Let $G$ be a reductive group that decomposes as a product
\eqn
G \cong \prod_{j \in J} {\rm{Res}}_{F_j/F}  \widetilde G_j
\label{eqn:dec-Gad}
\xeqn
for some finite index set $J$ and reductive groups $\widetilde G_j$ defined over a finite separable extension $F_j$ of $F$ for every $j \in J$. 
Recall that both $G^{\rm sc}$ and $G^{\rm ad}$ always have such a decomposition where the factors $\widetilde{G}_j$ are absolutely almost simple.
We recall the corresponding decomposition of the combinatorial data introduced above attached to $G$ following \cite[Section 4.2]{HainesRicharz:TestFunctionsWeil}, which discusses the case of Weil-restricted groups, but all the constructions are also commute with taking products, compare also \cite[Section 3]{HainesRicharz:Smoothness} and \cite[Section 3.5]{HainesRicharz:Normality}.

By the discussion in \cite[Section 4.2]{HainesRicharz:TestFunctionsWeil} we find 
\eqn
T \cong  \prod_{j \in J} {\rm{Res}}_{F_j/F} \widetilde T_j \qquad \text{and} \qquad B \cong  \prod_{j \in J} {\rm{Res}}_{F_j/F} \widetilde B_j
\label{eqn:dec-TBad}
\xeqn
under \eqref{eqn:dec-Gad}, where $\widetilde T_j$ and $\widetilde B_j$ are maximal $F_j$-tori and Borel subgroups of $\widetilde G_j$, respectively. Moreover, using \cite[Lemma 4.1]{HainesRicharz:TestFunctionsWeil} we have
\eqn
X_*(T)_I \cong \bigoplus_{j \in J} X_*(\widetilde T_j)_{I_j},
\label{eq:dec-cochars}
\xeqn
where $I_j = \Gal(F^{\rm{sep}}/F_j) \subseteq I$. 
The decomposition \eqref{eq:dec-cochars} restricts to an isomorphism of partially ordered monoids on dominant cocharacters.
Also, there is a bijective correspondence between maximal $F$-split tori $S \subseteq G$ and tuples $(\widetilde{S}_j)_{j \in J}$ of maximal $F_j$-split tori $\widetilde{S}_j \subseteq \widetilde{G}_j$. Namely, $S$ decomposes as 
$
	S = \prod_{j \in J} S_j
$
such that the base change $S_{j, F_j} = \widetilde{S}_j$, compare the discussion in \cite[Section 4.2]{HainesRicharz:TestFunctionsWeil}.
In a similar fashion, by \cite[Proposition 4.6]{HainesRicharz:TestFunctionsWeil} the (enlarged) Bruhat-Tits building decomposes (as a polysimplicial complex) as
\[
	\BTB(G,F) = \prod_{j \in J}  \BTB(\widetilde{G}_j, F_j)
\]
equivariantly for the action of $G(F) = \prod_{j \in J} \widetilde{G}_j(F_j)$ which moreover identifies apartments 
\[
	\BTA(G,S,F) = \prod_{j \in J} \BTA(\widetilde{G}_j, \widetilde{S}_j, F_j)
\]
equivariantly for the action of the Iwahori-Weyl groups
\eqn
	W(G,S) = \prod_{j \in J} W(\widetilde{G}_j, \widetilde{S}_j).
	\label{eq:dec-weyl}
\xeqn
In particular, facets $\bbf$ in $\BTA(G,S,F)$ correspond bijectively to tuples of facets $(\bbf_j)_{j \in J}$ where $\bbf_j$ is a facet in $\BTA(\widetilde{G}_j, \widetilde{S}_j, F_j)$. 
Then, the corresponding parahoric group schemes decompose accordingly as 
$
	\calG_{\bbf} \cong \prod_{j \in J} {\rm Res}_{\calO_j/\calO} \, \widetilde{\calG}_{j,\bbf_j},
$ 
compare \cite[Proposition 4.7]{HainesRicharz:TestFunctionsWeil}.
 
Let us record the following immediate consequence for future reference.
\lemm
\thlabel{lem:dec-pi1}
	Under the decomposition \eqref{eqn:dec-Gad} the algebraic fundamental group decomposes as
	\eqn
		\pi_1(G) \cong \bigoplus_{j \in J} \pi_1(\widetilde{G}_j)
		\label{eq:dec-pi1}
	\xeqn
	equivariantly for the $I$-action. 
	In particular, a prime $p$ divides 
	$\# \pi_1(G^{\der})$ if and only if
	$p \mid \# \pi_1(\widetilde G^{\der}_j)$ for some $j \in J$.
	Moreover,
	\eqn
	\pi_1(G)_I \cong \bigoplus_{j \in J} \pi_1(\widetilde {G}_j)_{I_j}.
	\label{eq:dec-pi1inv}
	\xeqn
\xlemm
\pf
	The decomposition \eqref{eq:dec-pi1} follows essentially by definition, \eqref{eq:dec-pi1inv} then follows from \eqref{eq:dec-cochars}.
	Moreover, the divisibility assertion then follows from the fact that $\pi_1({\rm Res}_{F_j/F} \widetilde G_j) \cong \pi_1(\widetilde{G}_j)^{[F_j \colon F]}$, so
	\[
	\# \pi_1(G^{\rm der}) = \prod_{j \in J} ( \# \pi_1(\widetilde{G}_j))^{[F_j \colon F]}.
	\]
\xpf

The above discussion also applies to the simply connected cover $G^{\rm{sc}}$ of an arbitrary reductive group $G$. Let now $J$ be the set of connected components of the relative Dynkin diagram of $G$. Then we have
\eqn
G^{\rm{sc}} \cong \prod_{j \in J_{\rm simple}} {\rm{Res}}_{F_j/F} \widetilde G_j^{\rm{sc}}
\label{eqn:dec-Gsc}
\xeqn
where $\widetilde G_j^{\rm{sc}}$ is absolutely almost simple over $F_j$ for every $j \in J$.
The groups $\widetilde G_j^{\rm sc}$ are called the absolutely almost simple factors of $G^{\rm sc}$ while the $G_j \defined {\rm Res}_{F_j/F} \widetilde G_j^{\rm sc}$ are the almost simple factors of $G^{\rm sc}$. Similarly, we have the absolutely simple factors $\widetilde G^{\rm ad}_j$ of $G^{\rm ad}$. 
In particular, we obtain an inclusion
\eqn
\bigoplus_{j \in J} X_*(\widetilde T_j^{\rm{sc}})_{I_j} \hookrightarrow X_*(T)_I.
\label{eqn:inclusion-Tsc-factors}
\xeqn

For every $j \in J$ let $\bar \mu_j^{\rm qm} \in X_*(T_j^{\rm sc})_{I_j}^+$ be the quasi-minuscule cocharacter.
\defi
The \emph{(factorwise) quasi-minuscule cocharacter} $\bar \mu^{\rm qm} \in X_*(T)_I^+$ is the image of $(\bar \mu_j^{\rm qm})_{j \in J} \in X_*(T^{\rm{sc}})^+_I$ under the inclusion \eqref{eqn:inclusion-Tsc-factors}.
\thlabel{def-fqm}
\xdefi

\subsection{Echellonage root system}
\label{sect--echellonage-roots}

Recall that  \cite{BruhatTits:Groups1} introduce a reduced root system $\Sigma$ inside $(X^*(T)^I)_{\bbR} \cong (X_*(T)_I)_{\bbR}$, called the \emph{\'echellonage root system}, such that its affine Weyl group is $W_{\rm aff}(\Sigma) = W_{\rm aff} = W(G^{\rm sc}, S^{\rm sc})$.
We denote its roots by $\Phi_{\Sigma}$, and its coroots by $\Phi_{\Sigma}^\vee$.
A set of simple \'echellonage coroots is given by the image of the simple absolute coroots under the projection ${\rm res}_I \colon X_*(T)\to X_*(T)_I$, while the \'echellonage roots can be constructed via a modified norm map, compare \cite[Section 3, Proposition 6.1]{Haines:Duality}.
In particular, our choice of basis $\Delta_{\rm abs}$ for the absolute root system yields a natural basis $\Delta_{\Sigma}$ for the \'echellonage root system. Moreover, $\Delta_{\Sigma}$ is in natural bijection with the set of simple relative roots $\Delta$, such that the simple \'echellonage and relative roots differ at most by a factor of $2^{\pm 1}$. 

As $X_*(T^{\rm sc})_I$ and $X_*(T^{\rm ad})_I$ are torsion free, they embed into $X_*(T)_{I, \bbR}$ and we have natural inclusions
\eqn
Q^{\vee} \cong X_*(T^{\rm{sc}})_I \subseteq X_*(T^{\rm{ad}})_I \subseteq P^{\vee},
\xeqn
of lattices inside $X_*(T)_{I, \bbR}$,
where $Q^{\vee}$ and $P^{\vee}$ are the coroot and coweight lattice of the root system $\Sigma$, compare \cite[Lemma 15]{HainesRapoport:Parahoric}.
%
%

By \cite[Corollary 1.8]{Richarz:Schubert}, the Bruhat order restricted to $X_*(T)^+_I$ can be more explicitly described using that $\bar \la \leq \bar \mu$ if and only if $\bar \mu - \bar \la$ is a linear combination of simple \'echellonage coroots with non-negative integral coefficients.
An element $\bar \mu \in X_*(T)_I^+$ is called \emph{minuscule}, if $\langle \bar \mu, \alpha \rangle \in \{0, \pm 1\}$ for all \'echellonage roots $\alpha \in \Phi_\Sigma$. 
The minuscule elements in $X_*(T)^+_I$ are precisely the minimal elements with respect to the Bruhat order on $X_*(T)^+_I$.
%

\lemm
	The norm map ${\rm N}_I \colon X_*(T)_I \to X_*(T)^I$ restricts to a monotone map $X_*(T)^+_I \to  X_*(T)^{+, I}$. 
\thlabel{lem:norm-map-monotone}
\xlemm
\pf
	For the first claim let $\bar \mu \in X_*(T)_I$ and let $\alpha \in X^*(T)$ be a positive root.
	Then by construction
	$\lan \alpha, \mu' \ran \geq 0$ if and only if $\lan N'_I \alpha, \bar{\mu} \ran \geq 0$.
	Hence, $\bar \mu$ is dominant if and only if $\mu'$ is dominant.
	For the second statement we note that $\bar \la \leq \bar \mu$ if and only if $\bar \mu - \bar \la$ is the sum of simple \'echellonage coroots, but the norm of a simple \'echellonage coroot is a sum of simple absolute coroots.
\xpf

\lemm
	The map $X_*(T)_I \to X_*(T^{\rm{ad}})_I$ is injective when restricted to any connected component for the Bruhat-order on $X_*(T)^+_I$ with image a connected component in the Bruhat order on the dominant coweights $(P^\vee)^+$ in the \'echellonage root system.
\thlabel{lem:identify-conn-comp-split}
\xlemm
\pf
	By definition, the difference of elements in a connected component is a linear combination of \'echellonage coroots and hence the map is injective. 
	Moreover, a minuscule $\bar \mu \in X_*(T)_I^+$ maps to a minuscule dominant coweight in $P^\vee$ by definition.
\xpf

A dominant cocharacter $\bar \mu \in X_*(T^{\rm sc})^+_I$ can be written as a linear combination $\bar \mu = \sum_{\alpha^\vee \in \Delta^\vee_{\Sigma}} c_{\alpha^\vee} \alpha^\vee$ with non-negative integer coefficients. 
The set of simple coroots $\alpha^\vee \in \Delta^\vee_{\Sigma}$ that appear with positive coefficients in the sum is called the support  ${\rm supp}(\bar \mu) \subseteq \Delta^\vee_{\Sigma}$.
For a general $\bar \mu \in X_*(T)^+_I$ there is a unique minuscule $\bar \la \in X_*(T)^+_I$ such that $\bar \la \leq \bar \mu$. 
In particular, $\bar \mu - \bar \la \in X_*(T^{\rm sc})^+_I$.
By a slight abuse of notation we set ${\rm supp}(\bar \mu) = {\rm supp}(\bar \mu - \bar \la)$.

By the classification of standard parabolic subgroups we get for any subset of the \'echellonage root basis $\Psi \subseteq \Delta_{\Sigma}$ a standard Levi subgroup $M_\Psi$ of $G$ with simple (relative) roots precisely given by $\Psi$. 
In particular, we get a Levi subgroup $M_{\bar \mu}$ from ${\rm supp}(\bar \mu) \subseteq \Delta_{\Sigma}$. 
Its set of absolute roots $\Phi_{\rm{abs}}(M_{\bar \mu},T)$ is then given by the preimage of $\Phi(M_{\bar \mu}, S)$ under the restriction map $\Phi_{\rm abs}(G, T) \to \Phi(G, S)$.
Equivalently, a basis of simple absolute coroots is given by the set of $\alpha^{\vee} \in \Delta^\vee_{\rm abs}$ such that its image $\bar \alpha^{\vee} \in \Delta^{\vee}_{\Sigma}$ appears in the expansion of $\bar \mu - \bar \la$ with positive coefficient.

In general, for $\mu \in X_*(T)^+$ the support of its projection $\bar \mu \in X_*(T)^+_I$ can be bigger than (the image of) the support of $\mu$, for example for ramified odd unitary groups the image of the non-zero minuscule coweights is the quasi-minuscule cocharacter.
Nevertheless, we still get an injection on fundamental groups.
\lemm
	As subsets of the set of $(G,T)$-absolute roots $\Phi(M_{\mu}, T) \subseteq \Phi(M_{\bar \mu}, T)$, 
	and the inclusion induces an injective group homomorphism
	$$\pi_1(M_\mu^{\der}) \hookrightarrow \pi_1(M_{\bar \mu}^{\rm der})$$
	on fundamental groups.
	In particular, if $p \mid \pi_1(M_{ \mu}^{\rm der})$, then $p \mid \pi_1(M_{\bar \mu}^{\rm der})$ for a prime number $p$.
	\thlabel{lem:compare-pi1}
\xlemm
\pf
	By construction, for $\la \in X_*(T)^+$ the unique minuscule cocharacter such that $\la \leq \mu$ with $\mu -\la = \sum_{\alpha \in \Delta} c_{\alpha} \alpha^\vee$ we have for their projections $\bar \la \leq \bar \mu$, as $\bar \mu - \bar \la = \sum_{\alpha \in {\rm supp}(\mu)} c_{\alpha} \bar \alpha^{\vee}$. 
	Hence, 
	$$ {\rm supp}(\bar \mu) = \res_I({\rm supp}(\mu)) \cup {\rm supp}(\bar \la).$$
	We get an injective group homomorphism using \cite[Lemma 5.3]{BiekerRicharz:normality} by noting that $\Q \Phi^{\vee}(M_{\mu}) \cap \Z \Phi^{\vee}(M_{\bar \mu}) ) = \Z \Phi^{\vee}(M_{\mu})$.
\xpf


\section{Affine Grassmannians and Schubert varieties}
We continue to use the setup from before. In particular, $S \subseteq T \subseteq B  \subseteq G$ is a reductive group $G$ over $F$ together with a Borel subgroup $B$ and a maximally split maximal torus $T$ and a maximal split torus $S$ all defined over $F$.
Moreover, we fix a facet $\bbf$ in the apartment $\BTA(G,S,F)$ in the (enlarged) Bruhat-Tits building.
We recall the definitions and some well-known properties of Schubert varieties in (partial) affine flag varieties.

\subsection{Definitions}

For a $k$-algebra $R$ we denote by $R\pot{t}$, respectively $R\rpot{t}$ the ring of formal power series, respectively Laurent series in the formal variable $t$.
The {\it loop group} $LG$ (respectively, {\it positive loop group} $L^+\calG_{\bbf}$) is the group-valued functor on the category of $k$-algebras $R$ defined by $LG(R)=G(R\rpot{t})$ (respectively, $L^+\calG_{\bbf}(R)=\calG_{\bbf}(R\pot{t})$). 
Then, $L^+\calG_{\bbf}\subset LG$ is a subgroup functor. 
One has $L^+\calG_{\bbf}=\lim_{i\geq 0} \calG_{\bbf,i}$ where $\calG_{\bbf,i}(R)=\calG_{\bbf}(R[t]/(t^{i+1}))$ is the group of $i$-jets.
Each $\calG_{\bbf,i}$ is a smooth, affine $k$-group scheme with connected fibers.

The \emph{twisted partial affine flag variety} is defined as the \'etale quotient
\[
	\Fl_{G,{\bbf}}\defined LG/L^+\calG_{\bbf},
\]
which is representable by an ind-projective ind-scheme over $k$. 
In the case that $\bbf = x$ is a special vertex the twisted partial affine flag variety is also called the \emph{(twisted) affine Grassmannian} 
\[
	\Gr_{G,x} = LG/L^+\calG_{x} = \Fl_{G,x}.
\]
For a second facet $\bbf' \subset B(G,F)$ the twisted partial affine flag variety comes equipped with a left action by $L^+\calG_{\bbf'}$ by left multiplication.
For an element $w \in L^+\calG_{\bbf'}(k) \backslash LG(k)/ L^+\calG_{\bbf}(k) \cong W_{\bbf'} \backslash W / W_{\bbf}$
the \emph{Schubert variety}
\[
	\Sv_w = \Sv_w(\bbf', \bbf)
\]
is defined as the reduced $L^+\calG_{\bbf'}$-orbit closure of $\dot{w} \cdot e$ for any representative $\dot{w} \in LG(k)$ of $w$. 
Recall when $\bbf = x$ there is a more explicit lift for $\bar \mu \in X_*(T) \cong W/W_x$ given by the image $\bar \mu(t)$ of $\mu(t) \in T(F')$ under the norm map $T(F') \to T(F) \subseteq LG(k)$ for a lift $\mu \in X_*(T)$ of $\bar \mu$.  
When $\bbf = \bbf' =x$ we also write 
$
	\Gr_{G,x,\leq \bar \mu} = \Sv_{\bar{\mu}(t)}.
$

Moreover, when $G = T$ is a torus, the Kottwitz homomorphism $\kappa_T$ induces a bijection $\Gr_{T,x}(k) \xrightarrow{\cong} X_*(T)_I$. In other words, the underlying reduced ind-scheme is given by $(\Gr_{T,x})_{\mathrm{red}} \cong \underline{X_*(T)_I}$. Then, the point $\bar{\mu}(t) L^+\calG_x \in \Gr_{G,x}(k)$ is the image of $\bar \mu$ under 
$
	X_*(T)_I \xrightarrow{\cong} \Gr_{T,x}(k) \hookrightarrow \Gr_{G,x}(k).
$

\subsection{Basic geometric properties of Schubert varieties}

The Schubert variety $\Sv_w$ is a projective $k$-variety of dimension $\ell({}_{\bbf'}w^\bbf)$, where for an element $v \in W$ we write $v^\bbf$ for the (unique) element of minimal length in its right $W_\bbf$-coset and  ${}_{\bbf'}w^\bbf$ for the (unique) element of maximal length in $\{(vw)^\bbf \colon v \in W_{\bbf'}\}$, compare \cite[Lemma 1.6 and Propostition 2.8]{Richarz:Schubert}.
The underlying set of $\Sv_w$ is the disjoint union of $L^+\calG_{\bbf'}$-orbits for all $v \in W_{\bbf'} \backslash W / W_\bbf$ with $v \leq w$.

In particular, for a special vertex $x$ the Schubert variety $\Gr_{G,x,\leq \bar \mu}$ has dimension $\langle \bar \mu, 2\rho \rangle$, where $\rho$ is the halfsum of the positive \'echellonage roots, compare \cite[Corollary 2.10]{Richarz:Schubert}. We have 
$$\Gr_{G,x,\leq \bar \mu}=\coprod_{\bar \la \in X_*(T)^+_I: \bar\la \leq \bar\mu} \Gr_{G,x,\bar\la},$$ 
where more precisely $\Gr_{G,x,\bar \la}$ is defined as the (sheaf-theoretic) image of the orbit map $L^+\calG_x\to \Gr_{G,x}$, $g\mapsto g\cdot \bar \mu(t)$. It is representable by an open subscheme of $\Gr_{G,x,\leq \la}$.
The left-$L^+\calG_x$-action on each Schubert variety factors through the group of $i$-jets $L^+\calG_x\to \calG_{x,i}$ for some $i\geq 0$.

\subsubsection{Normality results}\label{sect--Schubert-schemes-normality}
The following theorem shows that most Schubert varieties are normal:

\theo[{\cite[Theorem 0.3]{PappasRapoport:Twisted}}]
\thlabel{theo:normal}
	Assume either ${\rm char}(k)=0$, or ${\rm char}(k)>0$ and ${\rm char}(k)\nmid \#\pi_1(G^\der)$.
	Then, the Schubert variety $\Sv_w(\bbf', \bbf)$ is normal for all $w \in W_{\bbf'} \backslash W / W_{\bbf}$.
\xtheo

On the other hand, almost all Schubert varieties are non-normal in the case ${\rm char}(k)\mid \#\pi_1(G^\der)$:

\theo[{\cite[Theorem 1.1]{HainesLourencoRicharz:Normality}}]
\thlabel{theo:not-normal}
	Assume $G$ is almost simple, semisimple, splits over a tamely ramified extension of $F$ and ${\rm char}(k)\mid \#\pi_1(G)$.
	Then, the Schubert variety $\Sv_w(\bbf', \bbf)$ is normal for only finitely many $w \in W_{\bbf'} \backslash W / W_{\bbf}$.
\xtheo

\subsubsection{Connected components of the twisted partial affine flag variety}

The \emph{Kottwitz homomorphism} induces a bijection
\eqn
\pi_0(\Fl_{G, \bbf}) \xrightarrow{\simeq} \pi_1(G)_I,
\label{eqn:identification-conn-comp}
\xeqn
compare \cite[Theorem 0.1]{PappasRapoport:Twisted}.
Let $w_0 \in \Omega_{\bba}$ and choose a lift $\dot w_0 \in LG(k)$. Then left-multiplication by $\dot w_0$ induces an isomorphism 
$$ \Fl_{G, \bbf} \xrightarrow{\dot{w}_0 \cdot } \Fl_{G, \bbf}$$
such that the induced map on the connected components under the identifications \eqref{eqn:identification-conn-comp} and \eqref{eq:pi1-Omega} is given by mulitplication by $w_0$. 
We use the following (well-known) fact which we state explicitly for future reference. 
\lemm[{cf. \cite[Lemma 3.2]{BiekerRicharz:normality}}]
\thlabel{lem:translates-conn-comp}
\begin{enumerate}
	\item 
	Assume that $w_0 \in \Omega_{\bba}$ stabilizes a facet $\bbf'$ in the closure of $\bba$ and let $X  \subset \Fl_{G,\bbf}$ be an $L^+\calG_{\bbf'}$-invariant closed subscheme.
	Then $\dot{w}_0(X)$ is independent of the choice of lift $\dot w_0$ and $L^+\calG_{\bbf'}$-invariant.	 
	Moreover, if $X = \Sv_w = \Sv_w(\bbf', \bbf)$, then $\dot{w}_0(\Sv_w) = \Sv_{w_0w}$.
	\item 
	Let $\bar \la \in \Omega_{\bba}\cap X_*(T)^+_I$.
	Then multiplication by $\bar{\la}(t)$ induces an isomorphism 
	\eqn
		\Gr_{G,  x, \leq \bar \mu} \xrightarrow{\cong} \Gr_{G,x,\leq \bar \la + \bar \mu}.
	\xeqn
	In particular, the isomorphism type of $\Gr_{G,x, \leq \bar \mu}$ only depends on the image $\bar \mu^{\rm{ad}}$ in the adjoint group. 
\end{enumerate}
\xlemm
\pf 
\begin{enumerate}
	\item As $w_0$ stabilizes the facet $\bbf$, we have $\dot w_0 L^+\calG_{\bbf} \dot w^{-1}_0 = L^+\calG_{\bbf}$.
	When $X$ is $L^+ \calG_{\bbf}$-invariant, then clearly $\dot{w}_0(X)$ is independent of the choice of lift and $L^+\calG_{\bbf}$-invariant.
	When $X=S_w$ is a Schubert variety, then so is $\dot w_0(X)$, as Schubert varieties are precisely the $L^+\calG_{\bbf}$-invariant closed subvarieties.
	But then the equality $\dot{w}_0 L^+\calG_{\bbf} w = L^+\calG \dot{w}_0 w$ proves the claim.
	
	\item This is an immediate consequence of the previous assertion together with  \thref{lem:central-cochars}.
\end{enumerate}

\xpf

\subsection{Comparison with the adjoint quotient and simply connected cover}
Let us recall the following properties of Schubert varieties for products, restrictions of scalars and central extensions.

\lemm[{\cite[Lemma 3.2]{HainesRicharz:Smoothness}, \cite[Lemma 3.9]{HainesRicharz:Normality}}]
\thlabel{lem:sv-product-dec}
	Let $G$ be a reductive group such that $G = \prod_{j \in J} {\rm Res}_{F_j/F} \widetilde{G}_j$ for finite separable extensions $F_j/F$ and reductive groups $\widetilde G_j$ over $F_j$.
	Under the decomposition $w = (\widetilde w_j)_{j \in J} \in W_{\bbf'}\backslash W/W_{\bbf} = \prod_{j \in J} W_{\bbf_j'} \backslash (\widetilde W_j) / W_{\bbf_j}$ the corresponding Schubert variety decomposes as
	\eqn
		\Sv_{w}(\bbf', \bbf) = \prod_{j \in J} \Sv_{\widetilde w_j}(\bbf_j', \bbf_j).
	\xeqn
	In particular, $\Sv_{w}(\bbf', \bbf)$ is normal if and only if $\Sv_{\widetilde w_j}(\bbf_j', \bbf_j)$ is normal for every $j \in J$.
\xlemm

Recall by \cite[Section 3.4]{HainesRicharz:Normality} a map $G' \to G$ be a map of reductive groups which induces an isomorphism on adjoint quotients yields natural maps of the corresponding Bruhat-Tits buildings that restricts to a map $\BTA(G',S',F) \to \BTA(G,S,F)$ on apartments. 
Moreover, the map on the apartment induces a bijection on facets an is compatible with the induced map $W' = W(G',S') \to W(G,S) = W$ on Iwahori-Weyl groups.
In a similar fashion we get a map $\Fl_{G',{\bbf}} \to \Fl_{G, {\bbf}}$ that restricts to a map on Schubert varieties
\eqn
	\Sv'_{w'}(\bbf', \bbf) \to \Sv_w(\bbf', \bbf).
	\label{eqn:map-schubert-var-central-ext}
\xeqn

\lemm[{\cite[Proposition 3.5]{HainesRicharz:Normality}, \cite[Proposition 2.1]{HainesLourencoRicharz:Normality}}]
	The map \eqref{eqn:map-schubert-var-central-ext} 
	 is a finite birational universal homeomorphism which induces an isomorphism on normalizations.
	 Moreover, it is an isomorphism when ${G'}^{\rm der} \to G^{\rm der}$ is \'etale.
	 \thlabel{lem:map-schubert-var-central-ext}
\xlemm

In particular, by applying the two previous lemmas to the simply connected cover,
 we obtain for every $w = (\widetilde{w}_j^{\rm sc})_{j \in J} \in {W}_{\bbf'}\backslash (W_{\rm aff})/W_{\bbf} = \prod_{j \in J} {W}_{\bbf_j'}\backslash(\widetilde W^{\rm sc}_j)/W_{\bbf_j} $ finite birational universal homeomorphisms
\begin{equation}
	\prod_{j \in J} \Sv^{\rm sc}_{\widetilde w_j^{\rm sc}}(\bbf_j', \bbf_j) \to \Sv_w(\bbf', \bbf),
\end{equation} 
where the $\widetilde G_j^{\mathrm{sc}}$ are the absolutely almost simple factors of $G^{\rm sc}$. We obtain a similar map to the Schubert variety of the adjoint quotient of $G$. 
We conclude the section by a sufficient critertion for the map to the adjoint Schubert variety to be an isomorphism.
\coro
\thlabel{cor:comparison-sv-adjoint}
Let $w \in {W}_{\bbf'} \backslash W /W_{\bbf}$ with image $w^{\rm ad} \in {W}_{\bbf'}\backslash( W^{\rm ad})/W_{\bbf}$.
Assume one of the following two conditions:
\begin{enumerate}
	\item The central isogeny $G^{\rm der} \to G^{\rm ad}$ is \'etale, e.g. if $G^{\rm der}$ is adjoint.
	\item The Schubert variety 
	$\Sv^{\rm ad}_{w^{\rm ad}}(\bbf', \bbf)$ 
	 is normal.
\end{enumerate}
Then the map 
$\Sv_{w}(\bbf', \bbf) \to \Sv^{\rm ad}_{w^{\rm ad}}(\bbf', \bbf)$
 is an isomorphism.
\xcoro
\section{Transversal slices and a twisted Levi lemma}\label{sect--transversal-slice}

The Levi lemma is a tool to study singularities along strata via certain transversal slices in Schubert varieties in affine Grassmannians by reduction to certain Levi subgroups.  
It is proved in the split case in \cite{MalkinOstrikVybornov:MinimalDegenerations} over $\bbC$, \cite[Section 2.3]{Zhu:Introduction} over $k$ and \cite[Section 4]{Lourenco:Grassmannian} over $\Z$, and was generalized in \cite
[Section 5]{BiekerRicharz:normality} to situations with non-normal Schubert varieties.
In the tamely ramified case, a version of the Levi lemma was proved in \cite[Corollary 3.9]{Richarz:Diplom} and \cite[Proof of Theorem 4.11]{BessonHong:SmoothLocus} (for simply connected groups over $\bbC$). 

The main result of this section is a twisted version of the results in \cite
[Section 5]{BiekerRicharz:normality} for groups that split over a tamely ramified extension of $F$ in arbitrary characteristic.
As in the split case in \emph{loc.~cit.} some additional care is needed due to the existence of non-normal Schubert varieties.

\subsection{Descriptions in terms of invariants}

Now assume that $G$ splits over a tamely ramified extension, and let $F'$ be  
a tamely ramified extension of $F$ that splits $G$.
Then $F'$ is a cyclic extension of $F$ of degree $e$ with ${\rm char} (k) \nmid e$. We fix a generator $\tau \in \Gal(F'/F)$.
Let us fix a uniformizer $u \in F'$ such that $u^e = t$. Then $\tau u = \xi u$ for a primitive $e$-th root of unity $\xi \in k$. 
We assume that our fixed special vertex $x$ becomes hyperspecial over $F'$ under the embedding $\BTA(G,S,F) \to \BTA(G_{F'}, T_{F'}, F')$.

Recall that for an absolutely almost simple group $G$ a vertex in $\BTA(G,S,F)$ is called \emph{absolutely special} if its image under the simplicial embedding $\BTA(G,S,F) \to \BTA(G_{F'}, T_{F'}, F')$ is special (or equivalently hyperspecial) for the minimal splitting field $F'$ (or equivalently for all splitting fields $F'$) of $G$, compare \cite[Definition 2.2.1, Remark 2.2.2, Remark 4.6.1]{PappasZhou:Smoothness}.
For a general reductive group $G$, a vertex in $\BTA(G,S,F)$ is absolutely special if its projection $\BTA(\widetilde{G}_j^{\rm ad}, S_j, F_j)$ is an absolutely special vertex for every absolutely simple factor $\widetilde{G}_j^{\rm ad}$ of the adjoint quotient $G^{\rm ad}$ of $G$.

Recall that as $G$ is quasi-split absolutely special vertices in $B(G,F)$ always exist and all absolutely special vertices are special by \cite[Lemma 5.2]{HainesRicharz:Smoothness}. If none of the absolutely simple factors $\widetilde{G}_j^{\rm ad}$ is an odd unitary group, every special vertex is absolutely special.
In the case of odd unitary groups the special but not absolutely special vertices become hyperspecial only after a totally ramified extension of degree 4.

Let $H$ be the split reductive group over $k$ such that $G_{F'} \simeq H_{F'}$.
Our choice of special vertex $x$ gives rise to a pinned isomorphism
\[
(H, B_H, T_H, X_H) \cong (G_{F'}, B_{F'}, T_{F'}, X).
\]

Then the automorphism ${\rm id} \otimes \tau$ on ${\rm Res}_{F'/F} G_{F'} = {\rm Res}_{F'/F} H_{F'}$ induces an automorphism $\sigma = \sigma_0 \otimes \tau$ on $H_{F'}$ where $\sigma_0$ is induced by an automorphism of the pinned root datum of $H$. 
Then $G \cong ({\rm{Res}}_{F'/F} H_{F'})^{\sigma}$ and the parahoric group scheme $\calG_x$ can then be written as 
\[
\calG_{x} = {\rm Res}_{\calO'/\calO} (H_{\calO'})^{\sigma, o},
\]
where $(\str)^o$ denotes the (fiberwise) neutral component.
On the level of loop groups we obtain
\eqn\label{eq-loop-group-invariants}
LG = (LH)^{\sigma} \quad \text{and} \quad L^+\calG_x = (L^+H)^{\sigma, o}.
\xeqn
More generally, 
under the identification of Bruhat-Tits buildings $\BTB(G, F) \xrightarrow{\cong} \BTB(H,F')^\sigma$ the alcove $\bba$ in $\BTB(G, F)$ identifies with the set of $\sigma$-fixed points of an alcove in $\BTB(H, F')$ also denoted by $\bba$, similarly for facets $\bbf$ in the closure of $\bba$.
Then $\calG_{\bbf} = (\Res_{\calO'/\calO} \calH_{\bbf})^{\sigma, o}$ and $L^+\calG_{\bbf} = (L^+\calH_{\bbf})^{\sigma, o}$. 
In particular, the inclusion of loop groups induces a map
\eqn
	\Fl_{G,\bbf} \rightarrow (\Fl_{H, \bbf})^\sigma
	\label{eqn:inclusion-invariants-Gr}
\xeqn
which is equivariant for the loop group action.

\lemm
\begin{enumerate}
	\item 
	On topological spaces, the restriction of the map \eqref{eqn:inclusion-invariants-Gr} to $\Gr_{T,x} \to (\Gr_{T_H})^\sigma$ is given by the norm map $X_*(T)_I \to X_*(T)^{I}$. 
	In particular, it is injective if and only if $X_*(T)_I$ is torsion-free.
	\item 
	The restriction of \eqref{eqn:inclusion-invariants-Gr} to each connected component of $\Fl_{G,\bbf}$ is an isomorphism onto a connected component of $(\Fl_{H, \bbf})^\sigma$.
\end{enumerate}
\thlabel{lem:Gr-invariants}
\xlemm
\begin{proof}
	\begin{enumerate}
		\item This follows directly from the fact that the diagram 
		\begin{center}
			\begin{tikzcd}
				T(F') \arrow[r, "N"] \arrow[d, "\kappa_{T_{F'}}"] & T(F) \arrow[r, "\cong"] \arrow[d, "\kappa_T"] & T(F')^I \arrow[d] \\
				X_*(T) \arrow[r, two heads, "{\rm res}_I"] & X_*(T)_I \arrow[r, "N_I"] & X_*(T)^I
			\end{tikzcd}
		\end{center}
		commutes by construction, compare \cite[Proposition 11.1.1]{KalethaPrasad}.
		\item It suffices to treat the neutral component. In this case, the argument of \cite{PappasRapoport:LoopGroups} (which treats the simply connected case) carries over.
	\end{enumerate}
\end{proof}

\subsection{Strictly negative loop groups}
\label{sect--strictly-neg-loop-group}
In order to define transversal slices we recall the construction of strictly negative loop groups for tamely ramified grops following \cite[Section 3.6]{dHL:Frobenius}.
Let $H$ be an arbitrary (split) reductive group over $k$ with a Borel pair $T_H \subseteq B_H \subseteq H$.
The \emph{negative loop group} $L^-H$ is given by the functor on $k$-algebras 
\[ 
	R \mapsto H(R[t^{-1}])
\] 
and the \emph{strictly negative loop group} is then defined as 
\[
	L^{--} H \defined \ker\left( L^{-1} H \xrightarrow{t^{-1} \mapsto 0} H \right).
\]

Let $Z_H$ be the kernel of the map $H^{\rm sc} \to H^{\der}$. We then have an exact sequence 
\[
	1 \rightarrow Z_H \rightarrow H^{\rm sc} \rightarrow H
\]
which induces an exact sequence on strictly negative loop group functors
\[
	1 \rightarrow L^{--}Z_H \rightarrow L^{--} H^{\rm sc} \rightarrow L^{--}H.
\]
More generally, let $B_H^{\rm opp} = T_H \ltimes U^{\rm opp}$ be the opposite (to $B_H$) Borel subgroup of $H$. 
Then for an alcove $\bba \subset \BTA(H, T_H, F)$ we set $L^{--} \calH_{\bba}  \defined L^{--}H \rtimes U^{\rm opp}$, and for a facet $\bbf$ in the closure of $\bba$ we define
\[
	L^{--} \calH_{\bbf} \defined \bigcap_{w \in W_{H, \bbf}} {}^wL^{--}\calH_{\bba}.
\]

Let us now come back to the setting considered before and let $H$ be the split form of $G$ defined over $k$. 
For a facet $\bbf$ in the closure of $\bba$ the strictly negative loop group for $\calG_{\bbf}$ is then defined as
\eqn\label{defn-twisted-negative-loop-group}
L^{--} \calG_{\bbf} \defined (L^{--}\calH_{\bbf})^{\sigma, o}.
\xeqn
Note that the construction does not depend on the choice of extension $F'/F$ as above.
\lemm
	The strictly negative loop group is representable by an ind-affine group ind-scheme.
	The composition $L^{--}\calG_{\bbf}\subset LG \to \Fl_{\calG,\bbf}$ is representable by a quasi-compact open immersion.
 \thlabel{lem:rep-strictly-neg-loop-group}
\xlemm
\pf
	The first assertion is clear by construction.
	We note that moreover by construction $L^{--} \calG_\bbf$ is connected and hence the map factors through the neutral component.
	The second claim now follows as in the absolutely almost simple case in \cite[Corollary 3.9]{HainesLourencoRicharz:Normality} 
\xpf

\subsection{Transversal slices at special level}
We continue to use the setup from above. In particular, let $F'$ be a tamely ramified extension that splits $G$ such that the special vertex $x$ becomes hyperspecial over $F'$.
In particular, when $x$ is absolutely special we can take $F'$ to be the minimal splitting field of $G$, but when $G$ is an odd unitary group and $x$ is a special but not absolutely special vertex we instead take a totally ramified degree 4 extension $F'/F$.   

\defi\thlabel{transversal-slice-defi}
For $\bar \la, \bar \mu\in X_*(T)_I^+$ the \textit{transversal slice of $\Gr_{G,x,\leq \bar \mu}$ at $\bar \la(t)$} is the étale subsheaf of $\Gr_{G,x}$ defined as
\eqn\label{transversal-slice-defi:eq}
	\Gr_{G,x,\leq \bar \mu}^{\bar \la}\defined L^{--}\calG_x \cdot\bar \la(t)\cap \Gr_{G,x,\leq \bar \mu}.
\xeqn
\xdefi
We denote by $\la', \mu' \in X_*(T)^{I,+}$ the norms of $\bar \la$ and $\bar \mu$, respectively.
Then using \thref{lem:Gr-invariants} and \thref{lem:rep-strictly-neg-loop-group} we obtain an open and closed immersion
\eqn\label{transversal-slice-invariants:eq}
\Gr_{G,x,\leq \bar \mu}^{\bar \la} \hookrightarrow (\Gr_{H, \leq \mu'}^{\la'})^{\sigma}
\xeqn
and its image (if non-empty) is precisely a connected component.
The following lemma is the twisted version of \cite[Proposition 4.2]{BiekerRicharz:normality}, a similar statement is shown in \cite[Corollary 3.9]{Richarz:Diplom}. 

\prop\thlabel{transversal-slice-lemm}
For $\bar \la, \bar \mu\in X_*(T)_I^+$ the sheaf $\Gr_{G,x,\leq \bar \mu}^{\bar \la}$ is representable by a locally closed subscheme of $\Gr_{G,x}$. 
It is non-empty if and only if $\bar \la\leq \bar \mu$.
In this case, the following properties hold:
\begin{enumerate}
	\item\label{transversal-slice-lemm:smoothly-equivalent}
		The scheme $\Gr_{G,x,\leq \bar \mu}^{\bar \la}$ is smoothly equivalent to the open subvariety $\cup_{\bar \nu\in X_*(T)_I^+: \bar \la\leq \bar \nu\leq \bar \mu}\Gr_{G,x,\bar \nu}$ in $\Gr_{G,x,\leq\bar \mu}$.
		More precisely, for each $i\geq 0$ such that the $L^+\calG_x$-action on $\Gr_{G,x,\leq \bar \mu}$ factors through the group of $i$-jets $\calG_{x,i}$ the maps in the diagram
		\eqn\label{transversal-slice-lemm:smoothly-equivalent:eq1}
			\Gr_{G,x,\leq \bar \mu}^{\bar \la}\overset{\; {\rm pr}}{\Gets} \calG_{x,i}\x \Gr_{G,x,\leq \bar \mu}^{\bar \la} \overset{{\rm act}\;}{\lr} \cup_{\bar \nu\in X_*(T)_I^+:\bar \la\leq \bar \nu\leq \bar \mu}\Gr_{G,x, \bar \nu}
		\xeqn
		are smooth and surjective.
		Here ${\rm pr}\colon (g,x)\mapsto x$ and ${\rm act}\colon (g,x)\mapsto g\cdot x$ are the projection and the action, respectively. 
	\item\label{transversal-slice-lemm:basics}
		The scheme $\Gr_{G,x,\leq \bar \mu}^{\bar \la}$ is affine, integral, geometrically unibranch, regular in codimension $1$ and of dimension $\lan2\rho,\bar \mu-\bar \la\ran$.
\end{enumerate}
\xprop
\pf
The representability, the non-emptiness criterion and the surjectivity of the maps in \eqref{transversal-slice-lemm:smoothly-equivalent} follows as in the proof of \cite[Proposition 4.2]{BiekerRicharz:normality}. Note that by \thref{lem:norm-map-monotone} we have $\la' \leq \mu'$, so $(\Gr_{H, \leq \mu'}^{\la'})^{\sigma}$ is non-empty using \emph{loc.~cit.} The smoothness of the maps in \eqref{transversal-slice-lemm:smoothly-equivalent} follows from the split case in \emph{loc.~cit.} by taking invariants using \eqref{transversal-slice-invariants:eq} and \cite[Proposition 3.5]{Edixhoven:NeronModels}.

\eqref{transversal-slice-lemm:basics}:
The properties ``geometrically unibranch'', ``reduced'' and ``regular in codimension $1$'' are local in the smooth topology \StP{034D}, \StP{0DQ1}. Hence the corresponding statements for $\Gr_{G,x,\leq \bar \mu}^{\bar \la}$ follow from \eqref{transversal-slice-lemm:smoothly-equivalent} and \cite[Theorem~1.1]{HainesLourencoRicharz:Normality}.
It remains to show that $\Gr_{G,x,\leq \bar \mu}^{\bar \la}$ is irreducible of dimension $\lan2\rho,\bar \mu- \bar \la\ran$, but the corresponding argument in the split case in  \cite[Proposition 4.2]{BiekerRicharz:normality} carries over verbatim.
\xpf

\subsection{A twisted Levi lemma}\label{sect--levi-lemma}
We deduce a twisted version of the Levi lemma from \cite[Section 3.3]{MalkinOstrikVybornov:MinimalDegenerations} adapted to our purposes, compare also \cite[Corollary 3.9 (3)]{Richarz:Diplom} and \cite[Section 5]{BiekerRicharz:normality}.

Let $M\subset G$ be a Levi subgroup containing $T$. The special vertex $x \in \BTA(G,S,F)$ is again a special vertex $x \in \BTA(M,S,F)$ and we consequently obtain a special parahoric group scheme $\calM_x$ over $\calO$ with generic fibre $M$. 
The closed immersion on twisted affine Grassmannians $\Gr_{M,x}\to \Gr_{G,x}$ induces a locally closed immersion
\eqn\label{immersion-Levi-slice}
	\Gr_{M,x,\leq \bar \mu}^{\bar \la} \lr \Gr_{G,x,\leq \bar \mu}^{\bar \la}
\xeqn
for any $\bar \la, \bar \mu\in X_*(T)^+_I$ with $\bar \la\leq \bar \mu$.
Note that every $B$-dominant cocharacter in $X_*(T)_I$ is $B\cap M$-dominant, so the transversal slice for the Levi subgroup is well-defined, see \thref{transversal-slice-defi}.

\lemm\thlabel{levi-lemma}
If $\bar\mu-\bar\la$ lies in the coroot lattice of the \'echellonage root system of $M$, then the map \eqref{immersion-Levi-slice} is an isomorphism. 
\xlemm
\pf
The proof of \cite[Lemma 5.1]{BiekerRicharz:normality} carries over verbatim, namely the map is a closed immersion of integral schemes of the same dimension by \thref{transversal-slice-lemm} \eqref{transversal-slice-lemm:basics} and hence an isomorphism.
\xpf

\coro
	If $\bar\mu-\bar\la$ lies in the coroot lattice of the \'echellonage root system of $M$,
	the open subschemes $\sqcup_{\bar \la \leq \bar \nu \leq \bar \mu} \Gr_{M,x, \bar \nu}$ and $\sqcup_{\bar \la \leq \bar \nu \leq \bar \mu} \Gr_{G,x, \bar \nu}$ of the Schubert varieties $\Gr_{M,x, \leq \bar \mu}$ and $\Gr_{G,x, \leq \bar \mu}$, respectively, are smoothly equivalent.
	\thlabel{cor:levi-lemma}
\xcoro
\pf
	This is a combination of \thref{transversal-slice-lemm}\eqref{transversal-slice-lemm:smoothly-equivalent} and \thref{levi-lemma}.
\xpf

As an immediate consequence we also obtain one implication in the assertion of our main theorem for all special vertices.
Assume now that $\bar \la$ is minuscule and that $M = M_{\bar \mu}$ is the Levi corresponding to the support of $\bar \mu - \bar \la$. Then by definition $\bar \mu - \bar \la$ lies in the coroot lattice of the \'echellonage root system of $M$.
\coro
	The Schubert varieties $\Gr_{M_{\bar \mu},x, \leq \bar \mu}$ and $\Gr_{G,x, \leq \bar \mu}$ are smoothly equivalent.
	In particular, both are normal provided that ${\rm char}(k) \nmid \pi_1(M^{\der}_{\bar \mu})$. 
	\thlabel{cor:suff-crit-special}
\xcoro

\section{Normality of quasi-minuscule Schubert varieties for tamely ramified groups}
\label{sect--quasi-minuscule}
One source of explicit non-normal Schubert varieties are quasi-minuscule Schubert varieties in affine Grassmannians for absolutely almost simple groups at absolutely special level, compare {\cite[Corollary 6.2]{HainesLourencoRicharz:Normality}, which was extended in \cite[Proposition 1.2]{BiekerRicharz:normality}.
	We generalize the results to not necessarily absoultely almost simple groups (still subject to Condition \eqref{eqn:cond-G}).

	\subsection{A normality criterion in terms of tangent spaces}
	To show the (non-)normality of Schubert varieties below we adapt calculations of tangent spaces from \cite[Section 6]{HainesLourencoRicharz:Normality}.
	We recall some relevant material. 
	Note that in part of \emph{loc.~cit.} it is assumed that $G$ is almost simple, but the arguments go through in our setting for general reductive groups that split over a tamely ramified extension of $F$.
	We continue to use the setup from the previous section. In particular, we fix a tamely ramified extension $F'$ of $F$ that splits $G$ together with a generator $\tau \in \Gal(F'/F)$. Moreover, we denote by $H$ the split form of $G$ such that $G = (\Res_{F'/F} H_{F'})^\sigma$, where $\sigma = \sigma_0 \otimes \tau$ for an automorphism $\sigma_0$ of the Dynkin diagram of $H$.
	Let $T_H$ and $Z_H$ be as in \refsect{strictly-neg-loop-group}.
	The simply connected cover $H^{\rm sc}$ of $H$ is then the split form of $G^{\rm sc}$.
	
	Let us fix a facet $\bbf \in \BTB(G,F)$.
	Using that the strictly negative loop group is an open neighborhood of the base point, c.f. \thref{lem:rep-strictly-neg-loop-group}, 
	we have
	\[
		T_e \Gr_{G,\bbf} = T_e L^{--} \calG_\bbf.
	\]
	Moreover, under this identification the map $\Fl_{G^{\rm sc},\bbf} \to \Fl_{G,\bbf}$ induces map on tangent spaces
	\[
		T_e L^{--}\calG_\bbf^{\rm sc} \to T_e L^{--} \calG_\bbf 
	\]
	coming from the natural homomorphism $G^{\rm sc} \to G$.
	Then we set
	\[ 
		L^{--} Z \defined  (L^{--}Z_{H, F'})^{\sigma, o} = \ker(L^{--}\calG_\bbf^{\rm sc} \to L^{--} \calG_\bbf),
	\]
	where the second equality follows from the left exactness of taking invariants as in \cite[Lemma 5.11]{HainesLourencoRicharz:Normality}.
	Let $\frakz_H$ be the Lie algebra of $Z_H$. Then $\frakz_H \neq 0$ if and only if ${\rm char}(k) \mid \#\pi_1(H^\der)$.

	Under the decomposition of $G^{\rm sc}= \prod_{j \in J} {\rm Res}_{F_j/F} \widetilde{G}_j^{\rm sc}$ its split form decomposes as 
	\[
		H^{\rm sc} = \prod_{j \in J} \prod_{\Gal(F_j/F)} \widetilde{H}_j^{\rm sc},
	\] 
	where $\widetilde{H}_j^{\rm sc}$ is the split form of $\widetilde{G}^{\rm sc}_j$.  
	The strictly negative loop group of $G^{\rm sc}$ then decomposes as 
	\[
		L^{--} \calG_x^{\rm sc} = \prod_{j \in J} L^{--}\widetilde{\calG}_{j,x_j}^{\rm sc}
	\]
	where $x_j \in \BTB(\widetilde{G}_j^{\rm sc}, F_j)$ is the projection of $x$ to $\BTB(\widetilde{G}_j^{\rm sc}, F_j)$ under the identification $\BTB(G^{\rm sc},F) \cong \prod_{j \in J} \BTB(\widetilde{G}^{\rm sc}_j,F_j)$. By definition, $x_j$ is again absolutely special. 
	Then $\widetilde{\calG}_{j,x_j}^{\rm sc}$ is the corresponding parahoric group scheme for $\widetilde G^{\rm}_j$.
	
	\lemm[{cf. \cite[Example 5.8]{HainesLourencoRicharz:Normality}}]
	The tangent space of $L^{--} \calG^{\sc}_x$ decomposes as
	\eqn
		T_e L^{--} \calG^{\rm sc}_x = \bigoplus_{i \geq 1}(\mathfrak{h}^{\rm sc}[u^{-i}])^{\sigma} = \bigoplus_{i \geq 1} \bigoplus_{j \in J}(\widetilde{\mathfrak{h}}_j^{\rm sc}[u_j^{-i}])^{\sigma_j}.
		\label{eqn:dec-tangent-space}
	\xeqn
	\thlabel{lem:dec-tangent-space}
	\xlemm
	Here, we denote by $\widetilde{\frakh}^{\rm sc}_j$ the Lie algebra of $\widetilde{H}^{\rm sc}_j$ and by $\sigma_j$ the image of $\sigma^{[F_j \colon F]} \in \Gal(F'/F_j)$ under the projection to $\Gal(F_j'/F_j)$, where $F_j'$ is the minimal splitting field of $\widetilde{G}^{\rm sc}_j$ with uniformizer $u_j$. 
	Let us denote by $e_j \defined [F_j' \colon F_j]$ and $d_j \defined [F' \colon F_j']$.
	Then $u_j = u^{d_j}$.  
	\pf
		For the second equality in \eqref{eqn:dec-tangent-space} we use that by construction $\sigma$ acts componentwise (in $J$) on the decomposition 
		\eqn 
			\frakh^{\rm sc} \cong \bigoplus_{j \in J} \left( \bigoplus_{\Gal(F_j/F)} \widetilde{\frakh}_j^{\rm sc} \right).
		\label{eqn:dec-hsc}
		\xeqn 
		The $\sigma$-action on the $j$-th summand restricts to a $\sigma_j$-action on $\widetilde{\frakh}_j^{\rm sc}[u^{-1}]$ that identifies the invariants 
		$\left( \bigoplus_{\Gal(F_j/F)} \widetilde{\frakh}_j^{\rm sc}[u^{-i}] \right)^\sigma = \widetilde{\frakh}_j^{\rm sc}[u^{-i}]^{\sigma_j}$ for every $i\geq 1$.
		As the order of the operation of diagram automorphism $\sigma_{0,j}$ on $\widetilde{\frakh}_{j}^{\rm sc}$ is the degree $e_j$ of the extension $F_j'/F_j$
		the space of invariants $\widetilde{\frakh}_j^{\rm sc}[u^{-i}]^{\sigma_j}$ vanishes whenever $i$ is not a multiple of $d_j$.
		Moreover, otherwise we have
		\[
			\left( \widetilde{\frakh}_j^{\rm sc}[u^{-i \cdot d_j}]  \right)^{\sigma_j} = 	\left( \widetilde{\frakh}_j^{\rm sc}[u_j^{-i}]  \right)^{\sigma_j}.
		\]
	\xpf
	
	Let us also recall the following (mild generalization of a) normality criterion for Schubert varieties of \cite{HainesLourencoRicharz:Normality} in terms of their tangent spaces.
	\prop[{\cite[Proposition 2.1, Corollary 5.12]{HainesLourencoRicharz:Normality}}]
		\thlabel{prop:normality-tangent-spaces}
		Assume that $G$ splits over a tamely ramified extension of $F$. 
		Let $\bbf', \bbf$ be a pair of facets in the closure of $\bba$ and let $w \in W_{\bbf'} \backslash W_{\rm aff} /W_{\bbf}$.  
		The following are equivalent:
		\begin{enumerate}
			\item The Schubert variety 
			$S_w = S_w(\bbf',  \bbf)$ is normal.
			\item The map 
			$S^{\rm sc}_w \to S_w$ 
			is an isomorphism.
			\item The map on tangent spaces 
			$T_e S^{\rm sc}_w \to T_e S_w$ 			
			is injective.
			\item 
			$(T_e L^{--} Z) \cap (T_e S_w^{\rm sc}) = 0$ as subspaces of $L^{--} \calG_\bbf$.
		\end{enumerate}
	\xprop
	\pf
		The equivalence of the first three assertions is \cite[Proposition 2.1]{HainesLourencoRicharz:Normality}. 
		The equivalence with the last assertion is shown in \cite[Corollary 5.12]{HainesLourencoRicharz:Normality} for absolutely almost simple groups, its proof generalizes verbatim to our setting.
	\xpf 
%
	\subsection{The quasi-minuscule Schubert variety for absolutely special level}
	\label{sect--proof-prop-intro}
	Our next goal is to discuss the normality of the quasi-minuscule Schubert varieties in characteristic dividing the order of the fundamental group of $G^{\rm der}$. 
	For the rest of this section we assume that $x = \mathbf{0}$ is absolutely special.
	Recall that when $G$ is absolutely almost simple and splits over a tamely ramified extension, by \cite[Corollary 6.2]{HainesLourencoRicharz:Normality} the quasi-minuscule Schubert variety is not normal. 
	We want to extend this result to certain not necessarily absolutely almost simple groups. 
	
	However, there exist tamely ramified groups $G$ such that even in characteristic ${\rm char}(k) \mid \# \pi_1(G^{\rm der})$ the quasi-minuscule Schubert variety is normal.
	This shows that \thref{thm:normality-criterion} does not hold for arbitrary (even almost simple) tamely ramified groups. 
	\prop
		\thlabel{prop:res-sl2-normal}
		Assume that ${\rm char}(k) = 2$ and let $F'/F$ an extension of degree 3.
		For 
		\[
			G = ({\rm Res}_{F'/F} \SL_2)/\Bmu_2
		\] 
		the quasi-minuscule Schubert variety $\Gr_{G,\mathbf{0},\leq \mu^{\rm qm}}$ is normal.  
	\xprop
	\pf
		Let $Z_{H^{\rm ad}} \cong  \Bmu_2^3 \subseteq H^{\rm sc} \cong \SL_2^3$ be the center of $H^{\rm sc}$. 
		Then the Lie algebra $\frakz_H$ of $Z_H \cong \Bmu_2$ embeds diagonally $\frakz_H \cong k \hookrightarrow k^3 \cong \frakz_{H^{\rm ad}}$ with $\sigma_0$ acting via a cyclic permutation of the coordinates.
		In particular, $\sigma$ acts by multiplication with $\xi^{-1}$ on $k[u^{-1}]$, so $(\frakz_H[u^{-1}])^\sigma = 0$.
		In this case, we can moreover see directly (e.g.\@ by the calculation in \cite[Appendix B]{HainesLourencoRicharz:Normality}) that 
		\[
			T_e \Gr_{G^{\rm sc}, \mathbf{0},\leq \bar \mu^{\rm qm}} = T_e \Gr_{\SL_2, \leq \bar \mu^{\rm qm}} = \mathfrak{sl}_2[t^{-1}] = (\frakh[u^{-1}])^\sigma.
		\]
		Hence,  $\Gr_{G,\mathbf{0},\leq \mu^{\rm qm}}$ is normal by \thref{prop:normality-tangent-spaces}.
	\xpf
	
	\rema
		More generally, a similar argument shows that for an almost simple and simply connected split group $\widetilde{G}$ with ${\rm char}(k) \mid \pi_1(\widetilde{G}^{\rm ad})$ we find $(\frakz_H[u^{-1}])^{\sigma} = 0$. 
		If we were to know that $T_e \Gr_{G^{\rm sc},\mathbf{0}, \leq \mu^{\rm qm}} = (\frakh [u^{-1}])^\sigma$, as in the $\SL_2$-case above, we could deduce normality of $\Gr_{G,\mathbf{0}, \leq \mu^{\rm qm}}$. 
		However, in general the tangent spaces of Schubert varieties in positive characteristic seem to be difficult to control, compare the discussion in \cite[Remark 4.4]{HainesLourencoRicharz:Normality}.
	\xrema
	
	Nevertheless, we can show that the quasi-minuscule Schubert variety is not normal for a large class of reductive groups when ${\rm char}(k) \mid \# \pi_1(G^{\rm der})$ that split over a tamely ramified extension.
	
	\prop
	\thlabel{prop:intro}
	Assume that $G$ splits over a tamely ramified extension of $F$ such that all of its almost simple simply connected factors are absolutely almost simple.
	Let $\mathbf 0 \in \BTA(G,S,F)$ be an absolutely special vertex.
	Let $\bar \la \in X_*(T)^+_I$ be a minuscule cocharacter and let $\bar \mu^{\rm qm}\in X_*(T)^+_I$ be the factorwise quasi-minuscule cocharacter.
	
	If ${\rm char}(k) \mid \pi_1(G^{\rm der})$, then the Schubert variety $\Gr_{G,\mathbf 0, \leq \bar \la + \bar \mu^{\rm qm}}$ is non-normal.
	\xprop
	
%
%
%
	\pf
As in the proof of \cite[Proposition 1.2]{BiekerRicharz:normality} we may reduce to the case $\la = 0$. 
As $Z$ is not \'etale there is an index $j_0 \in J$ such that the image $Z_{j_0} \subset G_{j_0}^{\rm sc}$ of $Z$ under the projection $G^{\rm sc} \to G_{j_0}^{\rm sc}$ is not \'etale. 
We get a surjective $\sigma$-equivariant map on Lie algebras
\[
	\frakz_{H}[u^{-d_{j_0}}] \to \frakz_{{j_0}, H_{j_0}}[u^{-d_{j_0}}].
\]
In particular, as $\frakz_{{j_0}, H_{j_0}}[u^{-d_{j_0}}] = \frakz_{{j_0}, H_{j_0}}[u_{j_0}^{-1}]$ contains a non-zero $\sigma$-invariant vector by \cite[Corollary 6.2]{HainesLourencoRicharz:Normality}, so does $\frakz_{H}[u^{-d_{j_0}}]$.
Using the calculation in the proof of \thref{lem:dec-tangent-space}, we see that
\begin{align*}
	(\frakz_{H}[u^{-d_{j_0}}])^{\sigma}  & \subseteq \bigoplus_{j \in J} (\frakh_{j}[u^{-d_{j}}])^{\sigma} =  \bigoplus_{j \in J, e_j = e_{j_0}} (\frakh_{j}[u^{-d_{j}}])^{\sigma}  \\
	&=  \bigoplus_{j \in J,e_j = e_{j_0}}  (\frakh_{j}[u_j^{-1}])^{\sigma} \subseteq \bigoplus_{j \in J} T_e \Gr_{\widetilde{G}^{\rm sc}_j, \mathbf{0}, \leq \bar \mu_j^{\rm qm}} = T_e \Gr_{{G}^{\rm sc}, \mathbf{0}, \leq \bar \mu^{\rm qm}}.
\end{align*}
For the equality in the first line we use the fact that $e_j \in \{1,2,3\}$, and the inclusion in the second line is \cite[Proposition 6.1]{HainesLourencoRicharz:Normality}.
\xpf
	
\section{A normality criterion for Schubert varieties in affine Grassmannians}
\label{sect--normality-criterion}
We finish the proof of our main \thref{thm:normality-criterion} that gives a criterion for the (non-)normality of Schubert varieties in affine Grassmannians (at absolutely special level) in terms of the order of the fundamental group of a certain Levi subgroup. 

\subsection{A Levi subgroup}
As before, for a reductive group $G$ and $\bar \mu \in X_*(T)^+_I$, let 
$\bar \la \in X_*(T)^+_I$ be the unique minuscule cocharacter with $\bar \la \leq \bar \mu$. 
Let $M_{\bar \mu}$ be the standard Levi subgroup of $G$ corresponding to $\rm{supp}(\bar \mu - \bar \la)$ regarded as a subset of the set of simple (\'echellonage) coroots $\Delta_{\Sigma}^\vee$ for $G$.   
In particular, $\bar \la$ is also the minuscule cocharacter such that $\bar \la \leq \bar \mu$ in the Bruhat order for $M_{\bar \mu}$. 
\lemm\thlabel{constr-levi-cochar}
\begin{enumerate}
	\item 
	Let $\bar \mu \in X_*(T^{\rm sc})^+_I$. Then ${\rm supp}(\bar \mu) = \Delta_\Sigma$ if and only if $\bar \mu \geq \bar \mu^{\rm qm}$, where $\mu^{\rm qm} \in X_*(T)^{\rm qm}_I$ is the factorwise quasi-minuscule cocharacter. 
	\item 
	In the setup above, $M = M_{\bar \mu}$ is the unique standard Levi subgroup of $G$ such that $\bar \la + \bar \mu^{\rm{qm}}_M \leq \bar \mu$ in the Bruhat order on $X_*(T)^+_I$ for $M$, where $\bar \mu^{\rm{qm}}_M$ is the factorwise quasi-minuscule cocharacter for $M$.  
\end{enumerate}	
\xlemm 
\pf
\begin{enumerate}
	\item 
	This is clear in the almost simple case and follows in general using  \eqref{eq:dec-cochars}.
	\item 
	By the previous assertion the property $\bar \la + \bar \mu^{\rm qm}_M \leq \bar \mu$  is equivalent to $\Delta_M = \supp(\bar \mu - \bar \la)$. 
\end{enumerate}
\xpf

In order to apply the results of the previous section to $G$ and as well as its Levi subgroups we replace Condition \eqref{eqn:cond-G} by the stronger condition that 
\eqn
	G \simeq \prod_{j \in J} {\rm Res}_{F_j/F} \widetilde G_j \qquad  \text{s.t.} \quad \begin{array}{l}
		\text{every $\widetilde G_j$ splits over a tamely ramified extension of $F_j$,  } \\
		\text{and all of its almost simple simply connected} \\
		\text{factors are absolutely almost simple.}
	\end{array} 
	\label{eqn:cond-G-strong}
\xeqn
The condition is satisfied for all split groups as well as all simply connected or adjoint groups in charactersitic at least 5 or that split over a tamely ramified extension.

\lemm\thlabel{lem:levi-cond}
	Assume that $G$ satisfies \eqref{eqn:cond-G-strong}, and let $\bar \mu \in X_*(T)^+_I$. Then $M_{\bar \mu}$ also satisfies \eqref{eqn:cond-G-strong}.
\xlemm
\pf
	As the class of groups satisfying \eqref{eqn:cond-G-strong} is clearly closed under products and taking restrictions of scalars, and Levi subgroups of $G$ are of the form $M_{\bar \mu} \simeq \prod_{j \in J} {\rm Res}_{F_j/F} \widetilde{M}_{\bar \mu_j}$ it suffices to treat the case $G = \widetilde{G}_j$. 
	Here, $\bar \mu = (\mu_j)_{j \in J}$ under the decomposition $X_*(T)^+_I \cong \bigoplus_{j \in J} X_*(\widetilde T_j)^+_{I_j}$.
	
	By a similar argument, it moreover suffices to show that for absolutely almost simple groups $G$ all the almost simple simply connected factors of $M_{\bar \mu}$ are absolutely almost simple. 
	This is clear for all split groups. 
	For $G$ that split over a tamely ramified extension we have $M_{\bar \mu} = G$ except possibly in the case that $G$ is an even unitary group, i.e. $G^{\rm sc} \simeq {\rm SU}_{2n}$ and $\bar \mu^{\rm ad}$ is in the non-neutral component. 
	But in this case ${\rm supp}(\bar \mu) = \{2, \ldots, n\}$ and $M_{\bar \mu}^{\rm sc} \simeq {\rm SU}_{2(n-1)}$, compare Figure \eqref{figure-dom-cochars-Bn}. 
\xpf
\rema
	It is not true more generally that every standard Levi subgroup of a group $G$ satisfying \eqref{eqn:cond-G-strong} also satisfies the condition.
	Namely, when $G$ is an adjoint ramified triality the Levi subgroup corresponding to one of the two vertices in the relative Dynkin diagram of $G$ is isomorphic to $({\rm Res}_{F'/F} \SL_2)/\Bmu_2$, where $F'/F$ is the degree 3 extension that splits $G$.
\xrema

\subsection{The main theorem}
Using the Levi lemma for tamely ramified groups and the fact that factorwise quasi-minuscule Schubert varieties are non-normal under condition \eqref{eqn:cond-G-strong}, we deduce that we can read off the normality of the Schubert variety for $\bar \mu$ from the order of the fundamental group of $M^{\rm{der}}_{\bar \mu}$.

The result at absolutely special level implies a similar statement for more general level.
The $\bar{\mu}$-admissible set is defined as 
\[
	{\rm Adm}(G, \bar{\mu}) \defined \{ w \in W \colon w \leq w_0 \bar{\mu}(t) \text{ for some } w_0 \in W_{\mathbf{0}} \}.
\]
Moreover, for facets $\bbf, \bbf'$ in the closure of $\bba$ we set
$
	{}_{\bbf'} {\rm Adm}(G, \bar{\mu})_{\bbf} \defined W_{\bbf'}\backslash W_{\bbf'} {\rm Adm}(G, \bar{\mu}) W_{\bbf} / W_{\bbf}.
$

\theo
\thlabel{thm:main}
	\begin{enumerate}
		\item Assume that $G$ satisfies \eqref{eqn:cond-G}.
		Let $\bar \mu \in X_*(T)^+_I$ such that ${\rm char}(k) \nmid \pi_1(M_{\bar \mu}^{\rm der})$.  
		Suppose that $\bbf$ contains a special vertex $x$ in its closure. 
		Then $\Sv_w(\bbf', \bbf)$ is normal for every $w \in {}_{\bbf'}{\rm Adm}(G, \bar \mu)_{\bbf}$.
		\label{it:thm-normality-adm-locus}
		\item Moreover, when $G$ satisfies \eqref{eqn:cond-G-strong} and $\bbf = \bbf' = \mathbf{0}$ is an absolutely special vertex, $\Gr_{G, \mathbf{0}, \leq \bar \mu}$ is normal if and only if  
		$
			{\rm char}(k) \nmid \pi_1(M_{\bar \mu}^{\rm der}).
		$
		\label{it:thm-normality-nec}
	\end{enumerate}
\xtheo
In particular, we obtain \thref{thm:normality-criterion} from the introduction as an immediate corollary.
\pf
	Using \thref{lem:sv-product-dec} and \thref{lem:dec-pi1} we reduce to the tamely ramified case $G = \widetilde{G}_j$ for some $j \in J$.
  	By \thref{cor:suff-crit-special}, the Schubert variety $\Gr_{G,x,\leq \bar \mu}$ is normal in this case.
  	Part \eqref{it:thm-normality-adm-locus} now follows as in the proof of \cite[Proposition 9.1]{HainesLourencoRicharz:Normality}.
  	
  	To establish part \eqref{it:thm-normality-nec} by \thref{constr-levi-cochar} and \thref{cor:suff-crit-special} it suffices to show that $\Gr_{M_{\bar \mu}, {\mathbf 0}, \leq \bar \mu}$ is non-normal. 
  	But this follows from \thref{prop:intro} which is applicable due to \thref{lem:levi-cond}.
\xpf 
%
%

\subsection{A criterion for situations with only finitely many normal Schubert varieties}
We have the following criterion for groups such that each connected component of the partial affine flag variety only contains finitely many normal Schubert varieties.
For this let $G^{\rm sc} = \prod_{j} G_j^{\rm sc}$ be the decomposition of $G^{\rm sc}$ into almost simple factors, i.e. $G^{\rm sc}_j = {\rm Res}_{F_j/F} \widetilde{G}_j^{\rm sc}$ in the notation from above.
\prop
\thlabel{prop:criterion-finitely-normal}
Let $\bbf$ be a facet in $\BTA(G,S,F)$. The following assertions are equivalent: 
\begin{enumerate}
	\item 
	For every pair of facets $\bbf', \bbf$ each connected component of $\Fl_{G,\bbf}$ contains only finitely many normal $(\bbf', \bbf)$-Schubert varieties.
	\label{it:crit-finite-every-conn}
	\item 
	The twisted affine flag variety $\Fl_{G,{\bba}}$ contains a connected component which has only finitely many normal $(\bba, \bba)$-Schubert varieties.
	\label{it:crit-finite-some-conn}
	\item 
	For every $j \in J$ the kernel of $G^{\rm sc}_j \to G$ is not \'etale.
	\label{it:crit-finite-order-ker}
	\item 
	For every $j \in J$ the torsion part of the cokernel of the inclusion $X_*(T_j^{\rm{sc}}) \hookrightarrow X_*(T)$ has order divisible by ${\rm char}(k)$. 
	\label{it:crit-finite-order-tors-cochars}
	\item 
	For every $j \in J$ the characteristic ${\rm char}(k) \mid  \# \pi_1(M_j^{\rm{der}})$, where $M_j \subset G$ is the Levi subgroup corresponding to the set of relative roots of $G_j^{\rm sc}$.
	\label{it:crit-finite-order-fun-group}
\end{enumerate}
In particular, the affine flag variety $\Fl_{G, \bbf}$ contains only finitely many normal Schubert varieties if and only if the center of $G$ does not contain a non-trivial split torus and $G$ satisfies the above equivalent conditions.
\xprop
\pf
The last sentence is a consequence of the fact that $G$ not containing a non-trivial split central torus is equivalent to $\pi_0(\Fl_{G, \bbf}) = \pi_1(G)_I$ being finite.
Let us now show that \eqref{it:crit-finite-every-conn} - \eqref{it:crit-finite-order-fun-group} are equivalent.

\eqref{it:crit-finite-every-conn} $ \Rightarrow$ \eqref{it:crit-finite-some-conn}: This is clear.

\eqref{it:crit-finite-some-conn}  $ \Rightarrow$  \eqref{it:crit-finite-order-ker}: 
By \thref{lem:translates-conn-comp} it suffices to treat the neutral component of $\Fl_{G,{\bba}}$. 
By assumption, for every $j \in J$ there is $w_j \in \widetilde{W}_j^{\rm sc}$, where as above $\widetilde W^{\rm sc}_j = W(\widetilde{G}^{\rm sc}_j, \widetilde{S}^{\rm sc}_j)$, such that for $ w = (e, \ldots, e, w_j, e, \ldots, e) \in W_{\rm aff} = \prod_{j \in J} \widetilde{W}^{\rm sc}_j$ the Schubert variety $\Sv_w$ is nonnormal. 
Hence, $T_eL^{--} Z \cap T_e \Sv_{w} \neq 0$ by \thref{prop:normality-tangent-spaces}.
But this means in particular that 
$
	T_eL^{--} Z_j = {T}_e{L}^{--} Z \cap T_e L^{--} \calG_{\bba}  \supseteq T_eL^{--} Z \cap T_e \Sv_{w}
$ 
is nonzero. Thus, $Z_j$ cannot be \'etale.

\eqref{it:crit-finite-order-ker} $\Leftrightarrow$ \eqref{it:crit-finite-order-tors-cochars} $\Leftrightarrow$ \eqref{it:crit-finite-order-fun-group}:
As the kernel $Z_j$ is a finite group scheme of multiplicative type, it is \'etale if and only if the order of its geometric character group $X^*(Z_j)$ is not a multiple of ${\rm char}(k)$. 
By duality, the order of $X^*(Z_j)$ agrees with the order of the torsion part of the cokernel of $X_*(T_j^{\rm{sc}}) \hookrightarrow X_*(T)$.
Moreover, by \cite[Lemma 5.3]{BiekerRicharz:normality} the fundamental group $M_j^{\rm der}$ is precisely given by the torsion part of the cokernel of the inclusion $X_*(T_j^{\rm{sc}}) \hookrightarrow X_*(T)$.

\eqref{it:crit-finite-order-fun-group} $\Rightarrow$ \eqref{it:crit-finite-every-conn}:
As above, it suffices to treat the neutral component of $\Fl_{G,{\bbf}}$. 
We consider the central isogeny $\prod_{j \in J} M_j^{\rm der} \to G^{\rm der}$.
By applying \thref{lem:map-schubert-var-central-ext} we see that it suffices to show that the neutral component of $\prod_{j \in J} \Fl_{M_j^{\rm der}, \bbf_j}$ contains only finitely many normal Schubert varieties. But for every factor this follows from \cite{HainesLourencoRicharz:Normality}.

\xpf

\exam[The case $G = \SO_4$]
We consider the case $G = \SO_4 = (\SL_2 \times \SL_2)/\mu_2$ in characteristic 2. In this case, the derived subgroups of both Levi subgroups corresponding to the two almost simple factors are isomorphic to $\SL_2$.
We identify the dominant cocharacters of $\SO_4$ with $(2\Z_{\geq 0})^{2}  \sqcup (1 + 2\Z_{\geq 0})^{2}$.

By the above discussion the Schubert variety $\Gr_{G, \leq (a,b)}$ is normal if and only if at least one of the two components $a, b$ is either $0$ or $1$.
In particular, there are infinitely many normal Schubert varieties in both connected components.
\xexam

\section{The classification for absolutely almost simple tamely ramified groups}
\label{sect--classification}
In this section we finish the proof of the classification of normal Schubert varieties in twisted affine Grassmannians for (not necessarily semisimple) tamely ramified absolutely almost simple reductive groups $G$ in characteristic ${\rm char}(k)\mid \#\pi_1(G^{\rm der})$. 

Using \thref{thm:normality-criterion} our basic strategy is (as in the split case in \cite{BiekerRicharz:normality}) as follows:
\begin{enumerate}
	\item Determine all dominant cocharacters $\bar \mu \not \geq \bar \la + \bar \mu^{\rm qm}$, where $\bar \la \in X_*(T)^+_I$ is minuscule. For the cases relevant to us compare Appendix \ref{sec--bruhat-order}.
	\item For each such $\bar \mu \in X_*(T)^+_I$ determine if ${\rm char}(k) \mid \# \pi_1(M_{\bar \mu}^{\rm der})$, e.g. using \eqref{eqn:order-fundamental-group} or by determining the type of $M_{\bar \mu}$ and comparing its connection index with ${\rm char}(k)$.
	Note that it suffices to do the computation for the maximal normal and minimal non-normal Schubert varieties.
	This is done in the remainder of this section.
\end{enumerate}
%

\subsection{Proof of \thref{classification-normal-schubert-varieties} for semisimple groups}
As a first step we treat semisimple groups that split over a tamely ramified extension of $F$. Recall that \thref{classification-normal-schubert-varieties} for split semisimple groups is \cite[Theorem 1.1]{BiekerRicharz:normality}.
It remains to treat the case of non-split tamely ramified groups absolutely almost simple groups $G$.
Comparing with the list of \cite[§4]{Tits:Corvallis} such groups $G$ include the following cases:
\begin{itemize}
	\item even unitary groups, type $B\str C_n$ with $n \geq 3$, ${\rm char}(k) \neq 2$ and ${\rm char}(k)  \mid n$,
	\item odd unitary groups, type $C \str BC_n$ with $n \geq 1$ and ${\rm char}(k)  \mid 2n+1$,
	\item ramified $E_6$, type $F_4^I$ with ${\rm char}(k)  = 3$, and 
	\item ramified triality, type $G_2^I$ and ${\rm char}(k)  = 2$,
\end{itemize}
where the type refers to the type of \cite[§4]{Tits:Corvallis}. Note that the twisted orthogonal groups are excluded by our tame ramification hypothesis.
In all but the first case $\pi_0(\Gr_{G, \mathbf 0})$ is trivial by \cite[Lemma 4.2, Remark 4.3]{HainesRicharz:Smoothness}. In these cases \thref{classification-normal-schubert-varieties} reduces to \cite[Corollary 6.2]{HainesLourencoRicharz:Normality}. In particular, there are no normal Schubert varieties except for the trivial one.

It remains to consider the case of even unitary groups in the case ${\rm char}(k)  \mid 2n$ and  ${\rm char}(k)  \neq 2$,  in which case \thref{classification-normal-schubert-varieties} follows from the following result.
In this case the \'echellonage root system is of type $B_n$. 
We continue to use the notation from \cite{Bourbaki:Lie456}.

\prop
\thlabel{lem:even-unitary-case}
Let $G \simeq {\rm PU_{2n}}$ for $n \geq 3$ and let $\mathbf{0}$ be an absolutely special vertex.
Let $\bar \mu \in X_*(T)_I^+ \setminus \{0\}$. 
Then $\Gr_{G, \mathbf{0}, \leq \bar \mu}$ is normal if and only if  $\bar \mu = \cw{m}$ for some odd $1 \leq m \leq n$.
\xprop

\pf
By \thref{lem:cochars-Bn}, compare Figure \ref{figure-dom-cochars-Bn}, the dominant coweights $\cw{m}$ for $m$ odd are precisely the non-zero dominant coweights that are not at least $(\cw 1 +) \cw 2$. 
Hence, it suffices to show that $\Gr_{G, \mathbf 0, \leq \cw{m}}$ is normal. 
For $m \geq 3$ the support of $\cw{m} - \cw 1$ is $\{2, \ldots, n\}$, so $M_{\bar \mu}$ is of type $B\str C_{n-2}$. But as ${\rm char}(k)  \mid 2n$ and ${\rm char}(k)  \neq 2$, we find ${\rm char}(k) \nmid (2n-2)$. Hence,  $\Gr_{G, \mathbf{0}, \leq \cw{\tilde n}}$ is normal by \thref{thm:normality-criterion}.
\xpf
The case of not necessarily adjoint groups of type $B \str C_n$ now follows from \thref{cor:comparison-sv-adjoint}.

\subsection{The not necessarily semisimple case}
We now extend the results to the not necessarily semisimple case.
We continue to assume ${\rm{char}}(k) \mid \#\pi_1(G^{\der})$.
Whenever $G^{\der}$ is adjoint (e.g. in types $B_n$, $C_n$, $E_6$ and $E_7$) or when all Schubert varieties for the adjoint group for dominant cocharacters that are not at least of the form $\la + \mu^{\rm qm}$ for a minuscule cocharacter $\la$ are normal (as in types $A_n$, $C_n$, $E_6$, $E_7$ and $B \str C_n$), the theorem follows from the semisimple case together with \thref{cor:comparison-sv-adjoint}.
The remaining cases in type $D_n$  with ${\rm char}(k)  = 2$ are treated in the following.

\lemm
\begin{enumerate}
	\item 
	Assume that $n$ is even and $G^{\der} \not \simeq \SO_{2n, k}, \PSO_{2n, k}$. Then $\Gr_{G, \leq \mu}$ is normal for all $\mu$ with $\mu^{\rm ad} \leq \cw{n-1} + \cw{n}$.
	\item 
	Assume $\mu^{\rm ad} = \cw{n-1} + \cw 1$. Then $\Gr_{G, \leq \mu}$ is normal if and only if $G^\der \simeq \SO_{2n,k}$ for every $n \geq 4$, $G^\der \simeq \PSO_{2n,k}$ for $n$ odd, or $G^\der \not \simeq \SO_{2n, k}, \PSO_{2n, k}$ for $n\equiv 2\mod 4$.
	\item Assume $\mu^{\rm ad} = \cw{n} + \cw 1$. Then $\Gr_{G, \leq \mu}$ is normal if and only if $G^\der \simeq \SO_{2n,k}$ for every $n \geq 4$, $G^\der \simeq \PSO_{2n,k}$ for $n$ odd, or $G^\der \not \simeq \SO_{2n, k}, \PSO_{2n, k}$ for $4 \mid n $.
\end{enumerate}
\xlemm
By \thref{lem:cochars-Dn}, compare Figures \ref{figure-dom-cochars-D4}, \ref{figure-dom-cochars-Dn-even} and \ref{figure-dom-cochars-Dn-odd}, the coweights in the lemma are precisely those that are not at least of the form minuscule + quasi-minuscule.
Note that the first assertion of the lemma treats the connected components corresponding to minuscule $\la$ with $\la^{\rm ad} = \cw 1$, the corresponding cases for $G^{\rm der} \simeq \SO_{2n,k}, \PSO_{2n,k}$ are the first part of \cite[Lemma 7.6]{BiekerRicharz:normality}.
The second and third assertion treat the connected components for  minuscule $\la$ with $\la^{\rm ad} = \cw{n-1}$ and $\cw n$, respectively. 
In this case the corresponding assertions for $G \simeq \PSO_{2n,k}$ and for $\mu^{\rm ad} = \cw{n-1} + \cw 1$ for $G \not \simeq \SO_{2n, k}, \PSO_{2n, k}$ are contained in the second part of \emph{loc.~cit.}
\pf
\begin{enumerate}
	\item Assume that $n$ is even. Let $\cw 1 \lneq \mu^{\rm ad}$ and $\mu^{\rm ad} \leq  \cw{n-1} + \cw{n}$. 
	We have
	\begin{align*} 
		X_\ast(T \cap M^\der_{\mu})  = \Q  \Phi^\vee(M_{\mu}, T)\cap X_\ast(T) \simeq (0 \oplus \Q^{n-1}) \cap  \left( \Z^n_{\rm even} + (1/2, \ldots, 1/2) \right) = \Z^{n-1}_{\rm even}.
	\end{align*}
	Hence, $M_{\mu}$ has a simply connected derived subgroup in this case.
	\item In this case, $M_{\mu}$ has support $\{1, \ldots, n-1\}$, and $\Z \Phi^\vee(M_{\mu}, T) = \Z^n_{\Sigma = 0}$.
	It remains to treat $G \simeq \SO_{2n}$.
	We similarly have $X_*(T \cap M_\der) =  \Z^n_{\Sigma = 0}$, so $M_{\mu}^{\rm der}$ is simply connected.
	
	\item 	
	The case $G \simeq \SO_{2n,k}$ follows as in (2). 
	Similarly, the case $G \not \simeq \SO_{2n,k}, \PSO_{2n,k}$ follows as in \cite[Lemma 7.6]{BiekerRicharz:normality}.
\end{enumerate}
\xpf
This finishes the proof of \thref{classification-normal-schubert-varieties}.

\section{Towards a classification of normal Schubert varieties in the affine flag variety}
\label{sect--general-level}
The question which Schubert varieties are normal has an analogue also for Schubert varieties in (partial) affine flag varieties for other parahoric level structures. 
Note that in the case $G = \PGL_2$ the results of \cite[Corollary 6.7]{HainesLourencoRicharz:Normality} give a complete answer, but already for $G = \PGL_3$ the situation seems to be much more complicated, compare \cite[Section 6.4]{HainesLourencoRicharz:Normality}.

While our methods do not seem to be enough to give a full classification of normal Schubert varieties in the full affine flag variety, we do get new information about the normality of certain Schubert varieties.
In particular, we are able to classify all Iwahori-Schubert varieties in the affine Grassmannian in type $A_n$.
Additionally, we can classify all normal (Iwahori-)Schubert varieties in the affine flag variety for $\PU_3$. 


\subsection{Iwahori Schubert varieties in the affine Grassmannian of type $A_n$}

We classify normal Iwahori-Schubert varieties in the affine Grassmannian for a split almost simple group of type $A_n$.
Recall that by \thref{lem:translates-conn-comp} every Iwahori-Schubert variety is isomorphic to one in the neutral component via translation by a length 0 element in $W$.
Hence, it suffices to treat the neutral component of $\Gr_G$. 
Recall that Iwahori-orbits in the affine Grassmannian are parametrized by $W/W_0 \cong X_*(T)_I$.

As a first step, we consider $G\simeq \PGL_{n+1}$ and compute the translates of $n \cw 1$ and $n \cw n$ into the identity component in the coweight lattice.
Recall that the corresponding Schubert varieties in $\Gr_{\PGL_n}$ are maximal normal ones in the classification of \cite[Theorem 1.1]{BiekerRicharz:normality}, c.f. \thref{classification-normal-schubert-varieties}.

Let $w_0 \in \Omega_{\bba}$ the image of $\cw 1$ under the map 
\eqn
	X_*(T) \twoheadrightarrow \pi_1(G) \xrightarrow[\eqref{eq:pi1-Omega}]{\cong} \Omega_{\bba},
\xeqn
then $w_0$ is a generator of $\Omega_{\bba} \cong \pi_1(G) \cong \Z/(n+1)\Z$.
The element $w_0$ operates by left multiplication on $X_*(T) \xrightarrow[]{\cong} W/W_0$. 
\lemm
The translates of $n \cw{1}$ and $n \cw{n}$ to the neutral component are given by
\eqn
w_0 \cdot (n \cw 1) = (-1,\ldots,-1,n) \qquad \text{and} \qquad
w_0^{-1} \cdot( n \cw n )= (-n, 1, \ldots, 1),
\xeqn
where as usual $X_*(T^{\rm sc}) \cong \Z^{n+1}_{\Sigma = 0}$.
\xlemm
\pf
One checks that for $\mu = (\mu_0, \ldots, \mu_n) \in X_*(T)\cong \Z^{n+1}/(1,\ldots,1)$ we have
\eqn
w_0 \cdot \mu = (\mu_1, \ldots, \mu_n, \mu_0 + 1)
\xeqn
The assertion follows directly.
\xpf

The Bruhat order on $W$ descends to a partial order on $X_*(T)_I$ more explicitly given as follows.
\prop[Besson-Hong, {\cite[Proposition 6.2]{HainesLourencoRicharz:Normality}}]
	Let $\la, \mu \in X_*(T)_I$. Then $\la \leq \mu$ if and only if there exists a chain of $\mu_0, \ldots, \mu_r \in X_*(T)_I$ with $\mu = \mu_0$ and $\la = \mu_r$ such that there exists a positive root $\alpha_i$ such that $\mu_{i+1} = \mu_i - k  \alpha_i^\vee$ with $0 \leq k \leq \langle \alpha_i, \mu_{i} \rangle$  or such that $\mu_{i+1} = \mu_i - k \alpha_i^\vee$ with $0 \leq k < -\langle \alpha_i, \mu_i \rangle$.
	\thlabel{prop:Bruhat-order-cochars}
\xprop
Moreover, inside an orbit of the finite Weyl group on $X_*(T)_I$, the maximal element is the dominant one while the minimal element is the anti-dominant one.

As a next step we show that $(-1, \ldots, -1,n)$ and $(-n, 1, \ldots,1)$ are indeed the maximal elements in the set of cocharacters that are not as least as big as $\mu^{\rm qm}$.
\lemm
	\thlabel{lem:order-cochars-An}
	Let $\mu = (\mu_0, \ldots, \mu_n) \in X_*(T^{\rm sc}) \cong \Z^{n+1}_{\Sigma = 0}.$.
	Then $\mu$ satisfies at least one of the following assertions.
	\begin{enumerate}
		\item There are $0 \leq i,j \leq n$ with $\mu_i \leq -2$ and $\mu_j \geq 2$, in which case $\mu \geq \mu^{\rm qm}$.
		\label{it-22}
		\item For every $0 \leq i \leq n$ the components $\mu_i \leq 1$ and 
		\label{it-1}
		\begin{enumerate}
			\item $\mu_0 \leq 0$ in which case $\mu \leq (-n,1,\ldots,1)$, or 
			\label{it-1-0}
			\item $\mu_0=1$ and $\mu_n = -1$, in which case $\mu \geq \mu^{\rm qm}$.
			\label{it-1-1}
			\item $\mu_0 = 1$ and there exists a $1 \leq i \leq n$ with $\mu_i \leq -2$ in which case $\mu \geq \mu^{\rm qm}$, 
			\label{it-1-2}
		\end{enumerate}
		 \item For every $0 \leq i \leq n$ the components $\mu_i \geq -1$ and 
		 \label{it--1}
		 \begin{enumerate}
		 	\item $\mu_n \geq 0$ in which case $\mu \leq (-1,-1,\ldots,-1, n)$, or
		 	\label{it--1-0} 
		 	\item $\mu_0 = -1$ and there exists a $1 \leq i \leq n$ with $\mu_i \geq 2$ in which case $\mu \geq \mu^{\rm qm}$.
		 	\label{it--1--2}
		 \end{enumerate}
	\end{enumerate}
\xlemm	
Of course, the conditions in the lemma are not mutually exclusive, as \eqref{it-1-0} and \eqref{it--1-0} can be satisfied simultaneously, similarly, \eqref{it-1-1} and \eqref{it-1-2} can be satisfied at the same time.
\begin{proof}
	We first note that every $\mu \in \bbZ^{n+1}_{\Sigma = 0}$ satisfies one of the properties. 
	Namely, when $\mu$ does not satisfy \eqref{it-22}, we may without loss of generality assume that $\mu_i \leq 1$, in which case either $\mu_0 \leq 0$ or $\mu_0 = 1$. 
	We argue that the assertions in the lemma indeed cover all cases for $\mu = 1$.
	Namely, if \eqref{it-1-2} is not satisfies, then $|\mu_i| \leq 1$ for every $i \in I$. 
	We distinguish cases according to the value of $\mu_n$. Either $\mu_n = -1$ which is treated in \eqref{it-1-1}, or $\mu_n \geq 0$ which is treated in \eqref{it--1-0}.

	It remains to show the comparison to $\mu^{\rm qm}$, $(-n, 1, \ldots,1)$ or $(-1,\ldots,-1,n)$, respectively.
	It suffices to treat the cases \eqref{it-22} and \eqref{it-1} as \eqref{it--1} follows similarly.
	\begin{enumerate}
		\item As a first step using that $\mu \geq \mu^{\rm anti} = (m_0, \ldots, m_n)$ with $-m_0, m_n \geq 2$. Without loss of generality assume $m_n \geq -m_0$. If $m_n > -m_0$, then $n \geq 2$ and $m_1 < 0$. Let $\alpha = (0,1,0, \ldots,-1)$. So, $\langle \alpha, \mu^{\rm anti} \rangle = m_1 - m_n < m_n < -m_0-m_n$. 
		Hence, from \thref{prop:Bruhat-order-cochars} we obtain 
		$$\mu \geq (m_0, m_1 + (m_n+m_0), m_2, \ldots, m_{n-1}, -m_0) \geq (m_0, 0, \ldots,0, -m_0) \geq \mu^{\rm qm}.$$
		Where in the last step we use $\langle(1, 0, \ldots, 0,-1), (m_0, 0, \ldots, 0, -m_0)\rangle = 2m_0$. 
		\item 
		\begin{enumerate}
			\item We proceed iteratively for every $1 \leq i \leq n$ to construct $\mu^{(i)}$ with $\mu^{(0)} = (-n, 1, \ldots,1)$ such that in the $i$-th step we only modify the $0$-th and $i$-th entries and such that the $\mu^{(i-1)} \geq \mu^{(i)}$ and $\mu^{(i)}_i = \mu_i$, and such that $\mu^{(n)} = \mu$.  
			As the sum of all components is 0, we get 
			\eqn 
			n + \mu_0 =  \sum_{1 \leq i \leq n} 1 - \mu_i, 
			\label{eqn:sum-elements-mu}
			\xeqn
			where clearly a summand is nonzero if and only if $\mu_i \leq 0$.  
			Now consider the positive root $\alpha = (1, 0, \ldots, 0, -1, 0, \ldots, 0)$. Then  $\langle \alpha, \mu^{(i-1)} \rangle = \mu^{(i-1)}_0 - 1 > - (1 - \mu_i)$. Where the last inequality follows from \eqref{eqn:sum-elements-mu}.
			Then set $\mu^{(i)} = \mu^{(i-1)} + (\mu_i -1) \alpha$. 
			\item This is clear.
			\item We distinguish cases depending on the value of $\mu_n$. The case $\mu_n = -1$ was treated in \eqref{it-1-1}. We reduce the other two cases $\mu_n \geq 0$ and $\mu_n \leq -2$ to the case $\mu_n = -1$.  
			Let us first consider the case $\mu_n \geq 0$.  
			Then there exists a $\mu_i \leq -2$.
			Let $\alpha = (0, \ldots,0,1,0, \ldots, 0,-1)$ with the 1 in entry $i$. 
			Then $\langle \alpha, \mu \rangle = \mu_i - \mu_n \leq -( 2 + \mu_n)$. Hence, $\mu \geq \mu + (\mu_n + 1) \alpha$. But the latter has $-1$ as its last component.
			
			Let us now consider the case $\mu_n \leq -2$. Then there exists an index $1 \leq i \leq n-1$ with $\mu_i = 1$. Hence with $\alpha$ as before we find $\langle \alpha, \mu \rangle = 1 - \mu_n$, so $\mu \geq \mu - (\mu_n -1) \alpha$. 
		\end{enumerate}
	\end{enumerate}
\end{proof}

As a consequence of this discussion we obtain that all Iwahori-Schubert varieties in the affine Grassmannian for an absolutely simple group of type $A_n$ that do not contain the quasi-minuscule one are normal.
\coro
	Let $G$ be a (split) almost simple group of type $A_n$ and assume that ${\rm char}(k) \mid \# \pi_1(G^{\rm der})$. 
	For every $\mu \in X_*(T^{\rm sc})$ precisely one of the following is true:
	\begin{enumerate}
		\item $\mu \leq (-n, 1, \ldots, 1)$ or $\mu \leq (-1,\ldots,-1, n)$ in which case $\Sv_{\mu}(\bba, \mathbf 0)$ is normal.
		\label{it:coro-An-normal}
		\item $\mu \geq \mu^{\rm qm}$, in which case $\Sv_{\mu}(\bba, \mathbf 0)$ is non-normal.
		\label{it:coro-An-nonnormal}
	\end{enumerate}
\xcoro
\pf 
	Using \thref{lem:order-cochars-An} every $\mu \in X_*(T^{\rm sc})$ satisfies precisely one of the two conditions \eqref{it:coro-An-normal} and \eqref{it:coro-An-nonnormal}.
	In the second case $S_{\mu}(\bba, \mathbf 0)$ is non-normal by \cite[Corollary 6.2]{HainesLourencoRicharz:Normality}.
	
	In the first case, the Schubert variety is normal by combining \cite[Theorem 1.1]{BiekerRicharz:normality} (compare \thref{classification-normal-schubert-varieties}), \thref{lem:translates-conn-comp} and \thref{lem:order-cochars-An}.
\xpf

\rema
	By taking preimages under the (smooth) map $\Fl_{G, \bba} \to \Gr_{G,x}$ the (non-)normality of Iwahori-Schubert varieties in $\Gr_G$ also gives new information on (non-)normality of Iwahori-Schubert varieties in the affine flag variety. 
	For example, in the case $G = \PGL_3$ we get that all at most 6 dimensional Schubert varieties in the affine flag variety are normal while all at least 9 dimensional ones are not normal. For the seven and eight dimensional Iwahori-Schubert varieties our methods do not seem to be sufficient to give a complete answer.
	In particular, in the neutral component we find at least 70 normal Schubert varieties. For 24 Schubert varieties our methods do not seem to be enough to determine their normality, and all other ones are non-normal.
	\thlabel{rem:classification-pgl3}
	
	For general $n \geq 1$ our results show that there are non-normal Schubert varieties of dimension $\frac 12 (n^2 + 5n)$ and normal Schubert varieties of dimension $\frac 12 (3n^2+n)$ in the full affine flag variety for $\PGL_{n+1}$ with ${\rm char}(k) \mid (n+1)$.
	Note that for $n \geq 3$ the bound for the maximal dimension of normal Schubert varieties is strictly bigger than the bound for the minimal dimension of non-normal Schubert varieties.
\xrema

\subsection{The case of ramified odd unitary groups}
We collect some consequences in the case of tamely ramified odd unitary groups. In this case, the \'echellonage root system is of type $C_n$. 

\subsubsection{Special but not absolutely special level}
Recall that in ramified odd unitary groups there are special vertices that are not absolutely special.
In this case the quasi-minuscule Schubert variety is even smooth, (in contrast to the absolutely special case) by a theorem of Richarz, compare \cite{Arzdorf}.

Still, the non-normality of the quasi-minuscule Schubert variety at absolutely special level gives an explicit bound on the normal Schubert variety at the special but not absolutely special vertex.
\lemm
\thlabel{lem:odd-unitary-Schubertvarieties}
	Let $G \simeq \PU_{2n+1}$ be a tamely ramified odd unitary group in characteristic ${\rm char}(k) \mid (2n+1)$.
	Let $x$ be the special but not absolutely special vertex in the closure of the facet $\frak a$. 
	Then let 
	$$ \bar \mu = \begin{cases}
		2 \cw 1 & n = 1\\
		\cw 1 + 2 \cw n & n \geq 2
	\end{cases}.$$
	Then $\Gr_{G, x, \leq \bar \mu}$ is non-normal.
\xlemm
\pf 
	Let $w \in W$ be the longest element in the double coset $W_{\mathbf{0}} \bar{\mu}^{\rm qm}(t) W_{\mathbf{0}}$, where $\mathbf 0$ is the absolutely special vertex in the closure of $\bba$.
	Then one checks that $\bar \mu$ as in the assertion is the image of $w$ under the projection $W \to W_x \backslash W / W_x \xrightarrow{\cong} X_*(T)^+_I$.
	Then $S_w \subset \Fl_{G, \bba}$ and thus  $\Gr_{G, x, \leq \bar \mu}$ are non-normal by \cite[Proposition 2.3]{HainesLourencoRicharz:Normality}.
\xpf
In particular, in the case $n = 1$ this gives a full classification of all normal Schubert varieties in $\Gr_{G,x}$, while already for $n = 2$ we do not get an answer for the two Schubert varieties for $\bar \mu = 2 \cw 1$ and $\bar \mu = \cw 2$.
For bigger $n \geq 3$ this result is far from giving a full classification, compare Figure \ref{figure-dom-cochars-Cn}.

\subsubsection{Iwahori-Schubert varieties in the affine flag variety for $\PU_3$}
As in the example $G = \PGL_2$ in \cite[Corollary 6.7]{HainesLourencoRicharz:Normality} we can classify all (non-)normal Schubert varieties in the affine flag variety in the case $G = \PU_3$ in characteristic 3 such that $G$ splits over a ramified quadratic extension of $F$.
%
In this case, the results of \cite{HainesLourencoRicharz:Normality} and \cite{Arzdorf} are enough to classify all Schubert varieties in the affine flag variety.
\prop
	Let $G \simeq {\rm PU}_3$ and assume ${\rm char}(k) = 3$.
	Let $w \in {W}$. Then the Schubert $\Sv_w$ in the affine flag variety $\Fl_{G, \bba}$ is normal if and only if $\ell(w)$ is at most 2 or $w = s_0 s_1 s_0$, in which case $S_w$ is smooth, where $s_0$ is the affine reflection.
	\thlabel{prop:classification-PU3}
\xprop
\pf
There are precisely two 3-dimensional Schubert varieties in $\Fl_{G, \bba}$, one is the preimage of the quasi-minuscule Schubert variety at absolutely special level and hence non-normal, while the other is the preimage of the quasi-minuscule Schubert variety at special but not absolutely special level and hence smooth.
\xpf

\section{Normality of local models}
\label{sect--loc-model}
\index{Rework the local model part}
We use our results above to deduce a normality criterion and give a classification of normal local models as constructed by Pappas--Zhu \cite{PappasZhu:Kottwitz} and Levin \cite{Levin:Weil} in the mixed characteristic setting as well as Zhu \cite{Zhu:Coherence} and Richarz \cite{Richarz:AffGrass} in the equal characteristic setting, compare also \cite{HainesRicharz:Normality}.

In all cases the local models we consider here are defined as orbit closures inside certain Beilinson-Drinfeld affine Grassmannians. 
They are known to be normal whenever the characteristic of the residue field does not divide $\# \pi_1(\bbG^{\rm der})$, but due to the existence of non-normal Schubert varieties the local models fail to be normal in general, compare \cite[Section 9]{HainesLourencoRicharz:Normality}.
Let us also mention that in order to obtain normal local models one can correct the construction using $z$-extensions as in \cite{HePappasRapoport:semistable} or take seminormalizations, compare \cite{FHLR:singularities}.

In this section we denote by $\bbG$ a reductive group over a non-archimedean local field $\bbF$ (either in characteristic 0 or positive characteristic), with ring of integers $\mathbb O$ and residue field $\kappa$. 
Let $\breve \bbF$ be the completion of the maximal unramified extension of $\bbF$.
Let $\bbA \subseteq \bbS \subseteq \bbT \subseteq \bbG$ be a maximal $\bbF$-split torus $\bbA$, and a maximal $\breve \bbF$-split torus $\bbS$ with centralizer $\bbT$. 
Note that as $\bbG_{\breve{\bbF}}$ is quasi-split, $\bbT_{\breve \bbF} \subseteq \bbG_{\breve \bbF}$ is a maximal torus. Let us moreover fix a Borel subgroup $\bbB \subseteq \bbG_{\breve \bbF}$. 
As before, for a facet $\bbf \subset \BTA(\bbG, \bbA, \bbF)$ in the Bruhat-Tits building of $\bbG$ over $\bbF$ denote by $\calG_{\bbf}$ the corresponding parahoric group scheme. 
Moreover, let $\{\mu\}$ be a conjugacy class of (geometric) cocharacters of $\bbG$. 
Then $(\bbG, \{\mu\}, \calG_{\bbf})$ is an LM-triple in the sense of \cite{HePappasRapoport:semistable},
and we have in many cases an associated local model 
$	\mathcal M^{\rm loc} = \mathcal{M}^{\rm loc}(\bbG, \{\mu\}, \calG_{\bbf})$,
over the ring of integers $\bbO_{\bbE}$ of the reflex field $\bbE$ of $\{\mu\}$.

In ${\rm char}(\bbF) = 0$ the local model is defined by 
Pappas-Zhu \cite{PappasZhu:Kottwitz} when $\bbG$ splits over a tamely ramified extension of $\bbF$, Levin \cite{Levin:Weil} for groups $\bbG \simeq {\rm Res}_{F_1/F} \bbG_1$ for a finite extension $\bbF_1/\bbF$ and a reductive group $\bbG_1$ over $\bbF_1$ that splits over a tamely ramified extension of $\bbF_1$, and in \cite{FHLR:singularities} for groups including all adjoint and quasi-split groups with tamely ramified absolutely simple factors.
In the case ${\rm char}(\bbF) > 0$ the local model $\calM^{\rm loc}$ is constructed by Zhu \cite{Zhu:Coherence} and Richarz \cite{Richarz:AffGrass}.
In all cases it is constructed as a global Schubert variety inside a certain Beilinson-Drinfeld affine Grassmannian (in particular in the setting of \cite{FHLR:singularities} by $\calM^{\rm loc}$ we mean the orbit closure denoted by $M_{\calG, \mu}$ in \emph{loc.~cit.} and not its seminormalization). 

\rema
	It does not seem to be clear in general how to compare the constructions of \cite{Levin:Weil} and \cite{FHLR:singularities} when $\{\mu\}$ is not minuscule, compare \cite[Remark 5.13]{FHLR:singularities}.
	By a slight abuse of notation when we write $\mathcal{M}^{\rm loc}(\bbG, \{\mu\}, \calG_{\bbf})$ for $\bbF$ in characteristic 0, we always implicitly assume that it is defined in at least one of the settings of \cite{Levin:Weil} and \cite{FHLR:singularities}, and when it is defined in both settings our results are also satisfied for both constructions.
\xrema

The generic fiber of the local model is by construction given by the Schubert variety $\Gr_{\bbG, \leq \mu}$ for the group $\bbG \otimes_{\bbF} \bbF \rpot{\varpi - t}$, where $\varpi$ is a uniformizer of $\bbF$.   
Its special fiber can be embedded inside the affine flag variety for a certain group $G$ defined over $F \defined k\rpot{t}$ in positive characteristic, where $k \defined \kappa^{\rm alg}$, constructed as follows:
When ${\rm char}(\bbF) > 0$ set $G \defined \bbG_{\breve{\bbF}}$ to be the base change.
In the case ${\rm char}(\bbF) = 0$ we define $G$ to be the (base change to $F$) of the corresponding groups as in \cite[Section 4.1.1]{HainesRicharz:Normality} or \cite[Proposition 2.8, Section 5.2]{FHLR:singularities}, which are constructed from a certain lift of a parahoric model to $\bbO \pot{t}$.
In either case, we impose the condition that $G$ satisfies \eqref{eqn:cond-G}.
 
From the $\bbS$ and $\bbT$ we obtain similarly a maximal torus $T$ and is a maximal $F$-split torus $S$ in $G$.
In particular, we have a natural identification $X_*(T) \cong X_*(\bbT)$ compatibly with the actions of the absolute Weyl groups. 
Moreover, to the facet $\bbf$ we can in a natural way associate a facet $\bbf^\flat \subset \BTA(G,S,F) \subset \BTB(G,F)$, compare \cite[Section 4.1.3]{PappasZhu:Kottwitz} and a parahoric model $\calG_{\bbf^\flat}$ over $\calO$.

There is an embedding of the geometric special fiber $\calM^{\rm loc} \otimes_{\bbO_E} k \hookrightarrow \Fl_{G,{\bbf^\flat}}$, and the reduced subscheme of its image is the \emph{admissible locus} which is defined as the union of the Schubert varieties indexed by the \emph{admissible subset} ${\rm Adm}(G, \{\mu\}, \calG_{\bbf^\flat}) \subset W_{\bbf^\flat} \backslash W / W_{\bbf^\flat}$. 
Then \cite[Theorem 2.1]{HainesRicharz:Normality} and \cite[Corollary 9.2]{HainesLourencoRicharz:Normality} show that $\calM^{\rm loc}$ is normal if and only if all Schubert varieties that are contained in its generic or the special fiber are normal.
In \cite[Remark 2.2 (ii)]{HainesRicharz:Normality} and \cite[Proposition 9.1]{HainesLourencoRicharz:Normality} they show that this is always satisfied whenever $\bar \mu \in X_*(T)^+_I$ is minuscule.
We use our results on the normality of Schubert varieties to generalize this criterion.

As a first step we give a criterion for the normality of the generic fiber of the local model. 
Note that by construction the generic fiber of the local model is the Schubert variety for an unramified group. 
In particular, we can apply the (split case) of our results above. 
When ${\rm char}(\bbF) >0$ let $M_\mu$ be the Levi subgroup of $\bbG_{\bbF^{\rm sep}}$ corresponding to the support of $\mu$ as in \refsect{echellonage-roots}.

\prop
The generic fiber $\calM^{\rm loc} \times_{\Spec(\calO)} \Spec(\bbF)$ of the local model is normal if and only if ${\rm char}(\bbF) = 0$ or ${\rm char} (\bbF) \nmid \# \pi_1(M^{\der}_{\mu})$.
\thlabel{prop:normality-generic-fiber}
\xprop
\pf
This follows directly from the identification of the generic fiber of $\calM^{\rm 	loc}$ in positive characteristic with $\Gr_{G, \leq \mu}$, compare \cite[Lemma 3.4]{Richarz:AffGrass}, together with \thref{thm:normality-criterion}.
\xpf

\subsection{A normality criterion for local models}
In the following we denote for absolutely almost simple groups $\bbG$ by $\omega_i^\vee$ the basis of \emph{absolute} fundamental coweights of $\bbG$ from \cite{Bourbaki:Lie456}, and by $\bar{\omega}_i$ the fundamental coweights of the \'echellonage root system.  
 
\theo \thlabel{theo:normality-loc-mod}
	Let $(\bbG, \{\mu\}, \calG_{\bbf})$ be an LM-triple, such that $G$ satisfies Condition \eqref{eqn:cond-G} and assume that $\bbf$ contains a special vertex in its closure.
	\begin{enumerate}
		\item If $\mathrm{char}(\kappa) \nmid \# \pi_1(M_{\bar \mu}^{\der})$, then the local model $\calM^{\rm loc}$ is normal.
		\label{it:normality-loc-mod-crit}
		\item 
		\label{it:normality-loc-mod-class}
		Assume that $\bbG$ is absolutely almost simple, splits over a tamely ramified extension of $\bbF$ and that ${\rm char}(\kappa) \mid \# \pi_1(\bbG^{\der})$.
		Assume moreover that one of the conditions
		\begin{enumerate}
			\item $\bbG$ splits over an unramified extension of $\bbF$ and $(G, \{\mu\})$ appears in the list in \thref{classification-normal-schubert-varieties}.
			\label{it:normality-loc-mod-class-split}
			\item $\bbG_{\breve \bbF}^{\rm ad} \simeq \PU_{2n, \breve \bbF}$ and $\mu = \cw{2i-1}$ for some $1 \leq i \leq m$. 
			\label{it:normality-loc-mod-class-even}
			\item $\bbG_{\breve \bbF}^{\rm ad} \simeq \PU_{2n+1, \breve \bbF}$, $\mu \in \{\cw 1, \cw {2n}\}$, and $\bbf$ contains a special but not absolutely special vertex in its closure.
			\label{it:normality-loc-mod-class-odd}
		\end{enumerate}
		is satisfied.
		Then $\calM^{\rm loc}$ is normal.
		\item If $\bbf = \mathbf 0$ is an absolutely special vertex and $G$ satisfies \eqref{eqn:cond-G-strong}, then in \eqref{it:normality-loc-mod-crit} the condition $\mathrm{char}(\kappa) \nmid \# \pi_1(M_{\bar \mu}^{\der})$ is also necessary for $\calM^{\rm loc}(\bbG, \{\mu\}, \calG_{\mathbf{0}})$ to be normal. Moreover, for absolutely almost simple groups, the pairs $(\bbG, \{\mu\})$ as in \eqref{it:normality-loc-mod-class-split} and \eqref{it:normality-loc-mod-class-even} are the only situations such that $\calM^{\rm loc}(\bbG, \{\mu\}, \calG_{\mathbf{0}})$ is normal in this case.
	\end{enumerate}
\xtheo
The first assertion is a simultaneous generalization of both assertions in \cite[Remark 2.2(ii)]{HainesRicharz:Normality}. Namely, if $\mathrm{char}(\kappa) \nmid \# \pi_1(\bbG^{\der})$, then it also does not divide the order of the fundamental group of the derived subgroup of any Levi of $G$. Also, when $\bar \mu$ is minuscule, then $\pi_1(M_{\bar \mu}^{\der}) = 0$. 
\pf
\begin{enumerate}
	\item 
	The normality of the generic fiber follows from \thref{lem:compare-pi1} and \thref{prop:normality-generic-fiber}. 
	By \thref{thm:main} also all Schubert varieties in the admissible locus are normal. 
	We conclude using \cite[Theorem 2.1]{HainesRicharz:Normality}.
	\item 
	We again check that all Schubert varieties in the admissible locus are normal.
	As in the proof of \thref{thm:main} it suffices to show that $\Gr_{G, x, \leq \bar \mu}$ is normal. 
	In the split case \eqref{it:normality-loc-mod-class-split} we clearly have $\mu = \bar \mu \in X_*(T) = X_*(T)_I$. 
	Moreover, the cases \eqref{it:normality-loc-mod-class-odd} are the cases of exotic smoothness, in which the special fiber of $\calM^{\rm loc}(\bbG, \{\mu\}, \calG_{\mathbf{0}})$ is given by the Schubert variety for the quasi-minuscule cocharacter, which is smooth in this case.
	
	For case \eqref{it:normality-loc-mod-class-even} it remains to determine which $\mu$ reduce to $\bar \mu$ that appear in \thref{classification-normal-schubert-varieties} in the case $G^{\rm ad} \simeq {\rm PU}_{2n,F}$.
	But the reduction map $X_*(T^{\rm ad}) \cong \bbZ^{2n}/(1, \ldots, 1) \to X_*(T)_I \cong \bbZ^n$ is given by 
	$$ (\mu_1, \ldots, \mu_n, \mu_{n+1}, \ldots, \mu_{2n}) \mapsto (\mu_1 - \mu_{2n}, \ldots, \mu_n - \mu_{n+1}).$$
	One can now check that the preimage of $\bar{\omega}^\vee_m$ for $1 \leq m \leq  n$ is indeed given by $\{ \omega_{m}^{\vee}, \omega_{2n - m}^{\vee}\}$.
	This shows the claim.
	
	\item 
	This follows from \thref{thm:main} together with the calculation in \eqref{it:normality-loc-mod-class-even}.
\end{enumerate}
\xpf

\subsection{Groups of semisimple rank 1}
Using the classification of normal Schubert varieties in the affine flag variety for tamely ramified groups of semisimple rank 1 of \cite[Corollary 6.7]{HainesLourencoRicharz:Normality} in the case $G \simeq \PGL_2$ in characteristic 2 and \thref{prop:classification-PU3} for $G \simeq \PU_3$ in characteristic 3 we can classify normal local models for arbitrary parahoric level for these groups.
Examples of non-normal Schubert varieties in both cases are given in \cite[Example 9.3 and 9.4]{HainesLourencoRicharz:Normality}.

\theo
\thlabel{thm:class-locmod-ss1}
Let $(\bbG, \{\mu\}, \calG_{\bbf})$ be an LM-triple such that $\bbG$ is adjoint of $\breve \bbF$-rank 1 and its absolutely simple factor splits over a tamely ramified extension.
Assume moreover ${\rm char}(\kappa) \mid \pi_1(\bbG)$. 
Then $\calM^{\rm loc}(\bbG, \{\mu\}, \calG_\bbf)$ is normal if and only if the LM-triple appears in the following list.
\begin{enumerate}
	\item $\bbG_{\breve \bbF} \simeq {\rm Res}_{\breve \bbF'/\breve \bbF} \PGL_{2, \breve{\bbF}'}$ (so ${\rm char}(\kappa) = 2$) and 
	\begin{enumerate}
		\item $\mu = (0, \ldots, 0, \widetilde{\mu}_j, 0, \ldots,0)$ with $\widetilde{\mu}_j \in X_*(\widetilde{T})^+$ minuscule (for arbitrary $\calG_\bbf$ and $\bbF$), or
		\item $\mu = (0, \ldots, 0, \widetilde{\mu}_j, 0, \ldots,0, \widetilde{\mu}_{j'}, 0, \ldots,0)$ with $\widetilde{\mu}_j, \widetilde{\mu}_{j'} \in X_*(\widetilde{T})^+$ minuscule and $\calG_\bbf= \calG_\bba$ is an Iwahori group scheme, or
		\item $\mu = (0, \ldots, 0, \widetilde{\mu}_j, 0, \ldots,0)$ with $\widetilde{\mu}_j \in X_*(\widetilde{T})^+$ quasi-minuscule, ${\rm char}(\bbF) = 0$ and $\calG_\bbf = \calG_\bba$ is an Iwahori group scheme.
	\end{enumerate}
	\item $\bbG_{\breve \bbF}  \simeq {\rm Res}_{\breve \bbF'/\breve \bbF} \PU_{3, \breve{\bbF}'}$ (so ${\rm char}(\kappa) = 3$) and
	\begin{enumerate}
		\item $\mu = 0$, or
		\item $\mu = (0, \ldots, 0, \widetilde{\mu}_j, 0, \ldots,0)$ with $\widetilde \mu_j \in \{\cw 1, \cw 2\}$ is minuscule and 
		\begin{itemize}
			\item $\bbf = x$ is a special but not absolutely special vertex or 
			\item $\bbf = \bba$ is an alcove.
		\end{itemize}
	\end{enumerate}
\end{enumerate}
\xtheo

\pf
	We first check when all Schubert varieties $S_w$ for $w \in {\rm Adm}(G, \{\mu\}, \calG_{\bbf^\flat})$ are normal. 
	For absolutely special vertices this is contained in \thref{theo:normality-loc-mod}.
	For the special but not absolutely special vertex in the unitary case using \thref{lem:odd-unitary-Schubertvarieties} the Schubert variety in the special fiber of $\calM^{\rm loc}$ is normal if the image of $\mu \in X_*(T)^+$ under the projection $X_*(T) \cong (\bbZ_{\Sigma = 0})^{[F' \colon F]} \to X_*(T)_I \cong \bbZ$, where the map in each factor is given by $\la = (a,b,c) \mapsto \bar \la = a-c$ is either $\bar \mu = 0$ or $\bar \mu = \bar \mu^{\rm qm} = 1$. But the dominant cocharacters $\mu$ having this property are precisely those that have at most one non-zero component which is then given by a $\mu \in \{\cw 1, \cw 2\}$.
	
	For the Iwahori case we note that the maximal Schubert varieties in the admissible locus are indexed by the two elements in $W_{\rm aff}$ of length the dimension of $\Gr_{G, \leq \mu}$. 
	But by the explicit classification of normal Schubert varieties in the corresponding affine flag variety these are normal if and only if their dimension is at most two, 
	in which case either $\mu = 0$ or $\mu = \mu^{\rm qm} = (1,-1)$ in the $\PGL_2$-case or $\mu \in \{\cw 1, \cw 2\}$ in the $\PU_3$-case. 
	
	It remains to check in all cases if the generic fiber is also normal when ${\rm char}(F) > 0$. By inspection of the cases above we see that this is true except for the case $\mu = \mu^{\rm qm}$ in the $\PGL_2$-case. 
\xpf

\rema
	In particular, in the case $\bbG = \PGL_2$ in equal characteristic 2, the local model for the quasi-minuscule cocharacter at Iwahori level fails to be normal even though all Schubert varieties that are expected to appear in the special fiber are normal. 
	This is in contrast to the situation for absolutely special level where (at least for groups that satisfy Condition \eqref{eqn:cond-G-strong}, which includes all split groups) the normality of the Schubert variety in the special fiber implies the normality of the generic fiber in the local model by \thref{theo:normality-loc-mod}.	
\xrema

\rema 
	For more general tamely ramified groups $\bbG$ that have semisimple $\breve F$-rank 1 the local model is normal if $(\bbG^{\rm ad}, \{\mu^{\rm ad}\})$ appears in the list of \thref{thm:class-locmod-ss1}.
	However, as the quasi-minuscule Schubert variety in the special fiber might be normal, compare \thref{prop:res-sl2-normal}, there can be more normal local models in general. 
	For example, when $\bbF'/\bbF$ is a totally ramified degree 3 extension in residue characteristic 2, $\bbG = ({\rm Res}_{\bbF'/\bbF}\SL_2)/\Bmu_2$ and $\mu = (1,1,0)$ under the identification of $X_*(T)$ with elements in $\Z^3$ whose sum is even, the local model for absolutely special level is normal. 
\xrema

\rema
	Using the partial results in the $\PGL_3$-case in \thref{rem:classification-pgl3} we also get partial results for the Pappas-Zhu local model for $G^{\rm ad} \simeq \PGL_3$ in residue characteristic 3. 
	As all Schubert varieties in the affine flag variety of dimension at most 6 are normal we obtain for example that the Pappas-Zhu local model for $\mu = \mu^{\rm qm} = (1,0,-1)$ are normal for all facets $\bbf$ that are not vertices. 
	Moreover, at Iwahori level all local models of relative dimension at most 6 are normal.
\xrema

\appendix

\section{Combinatorics of the Bruhat order on dominant coweights}
\label{sec--bruhat-order}

We work out explicitly certain properties of the Bruhat order on dominant coweights in types $B_n, C_n$ and $D_n$ that are used above.
We use the notation of the previous sections and in particular denote by $\cw 1, \ldots, \cw n$  the fundamental coweights as in  \cite[Tables]{Bourbaki:Lie456}.
In particular, we give in  Figures \ref{figure-dom-cochars-Bn} - \ref{figure-dom-cochars-Dn-odd} certain initial segmets of the Bruhat order that are relevant in the above, extending corresponding discussion in \cite[Section 7]{BiekerRicharz:normality}.
As in \emph{loc.~cit.} the labels on the arrows are the support of the corresponding minimal degeneration and the correctness of the diagrams can be checked using the Stembridge classification of minimal degenerations \cite{Stembridge:Dominant}, compare also \cite[Theorem 6.3]{BiekerRicharz:normality}.

In type $B_n$ for $n \geq 3$ we have the following. 
\lemm
Let $\mu \in P^{\vee, +}$ with $\cw 1 \leq \mu$. Then $\mu \not \geq \cw 1 + \cw 2$ if and only if $\mu \leq \cw{\tilde n}$, where $\tilde n$ is the largest integer odd not bigger than $n$, i.e. $\tilde n \in \{n-1, n\}$ depending on the parity of $n$.
\thlabel{lem:cochars-Bn}
\xlemm
\pf
We check Figure \ref{figure-dom-cochars-Bn} using the Stembridge classification as in \cite[Lemma 7.4]{BiekerRicharz:normality}.
\xpf

\begin{figure}
	\begin{tikzcd}
		0 \arrow{rr}	& &  \cw 2 \arrow{rr}	&&		\ldots \\
	\end{tikzcd}
	
	\begin{tikzcd}
		\cw 1 \arrow{rr}{\{2,\ldots,n\}}	&&  \cw 3 \arrow{rr}{\{4, \ldots, n\}}\arrow{d}{\{1,2\}}	&& \cw 5  \arrow{r}\arrow{d} & \dots \arrow{r} & \cw{\tilde n}\arrow{d} &  \\	
		&& \cw 1 + \cw 2 \arrow{rr}\arrow{d} && \cw 1 + \cw 4 \arrow{r}\arrow{d} & \dots \arrow{r} &  \cw 1 + \cw{\tilde n -2} \arrow{d} \arrow{r} & \cdots \\
		&& \vdots && \vdots && \vdots & 
	\end{tikzcd}
	\caption{The Bruhat partial order on dominant coweights in type $B_n$ for $n \geq 3$. }
	\label{figure-dom-cochars-Bn}
\end{figure}

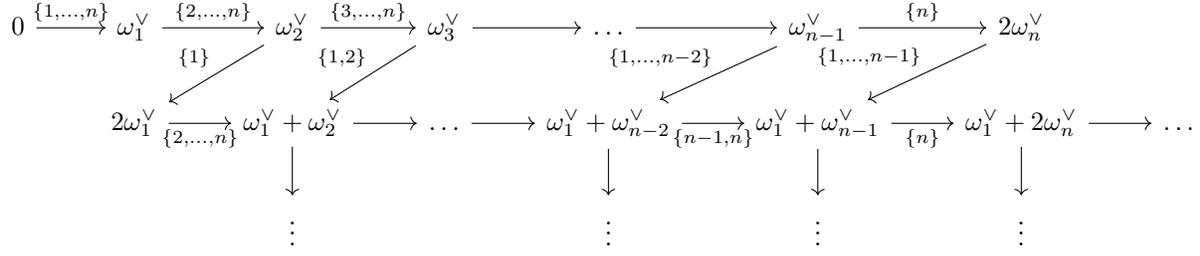
\begin{figure}
	\begin{tikzcd}
		0 \arrow{r}{\{1, \ldots, n\}}	&  \cw 1 \arrow{r}{\{2, \ldots, n\}}	& \cw 2  \arrow{r}{\{3, \ldots, n\}}\arrow{dl}[swap]{\{1\}} & \cw 3  \arrow{r}\arrow{dl}[swap]{\{1,2\}} & \dots \arrow{r} & \cw{n-1} \arrow{r}{\{n\}} \arrow{dl}[swap]{\{1, \ldots, n-2\}} &  2\cw{n}\arrow{dl}[swap]{\{1, \ldots, n-1\}}  &  & \\	
		&   2 \cw 1  \arrow{r}[swap]{\{2, \ldots, n\}}&  \cw 1 + \cw 2 \arrow{r}\arrow{d} & \dots \arrow{r} &  \cw 1 + \cw{ n -2} \arrow{d} \arrow{r}[swap]{\{n-1, n\}} & \cw 1 + \cw{n -1} \arrow{d}  \arrow{r}[swap]{\{n\}} &  \cw 1 + 2 \cw n  \arrow{r} \arrow{d}&  \ldots\\
		& & \vdots & & \vdots&   \vdots & \vdots &
	\end{tikzcd}
	
	\caption{The Bruhat partial order on dominant coweights in type $C_n$ for $n \geq 2$.
		The second connected component is obtained by shifting by $\cw n$.}
	\label{figure-dom-cochars-Cn}
\end{figure}

In type $D_n$ for $n \geq 4$ we have the following.
\lemm
\begin{enumerate}
	\item 
	The dominant coweight $\cw 1$ has a unique minimal degeneration $\cw 1 < \mu$ with $\mu = (1,1,1,0, \ldots,0)$, i.e. $\mu = \cw{3} + \cw{4}$ when $n = 4$ and $\mu = \cw 3$ when $n \geq 5$.
	\item
	When $n$ is even, a coweight $\mu \geq \cw 1$ satisfies $\mu \not \geq \cw 1 + \cw 2$ if and only if $\mu \leq \cw{n-1} + \cw n$.
	When $n$ is odd, a coweight $\mu \geq \cw 1$ satisfies $\mu \not \geq \cw 1 + \cw 2$ if and only if $\mu \leq 2 \cw{n-1}$ or $\mu \leq 2 \cw{n}$.
	\item 
	The connected components for the minuscule coweights $\la =\cw{n-1}$ or $\la = \cw{n}$ both contain a unique non-minuscule coweight which is not bigger than $\la + \cw 2$, given by $\cw 1 + \cw n$ and $\cw 1 + \cw{n-1}$, respectively.
\end{enumerate}
\thlabel{lem:cochars-Dn}
\xlemm
\pf
As before, we use the Stembridge classification to check that Figures \ref{figure-dom-cochars-D4} - \ref{figure-dom-cochars-Dn-odd} are correct.
\xpf

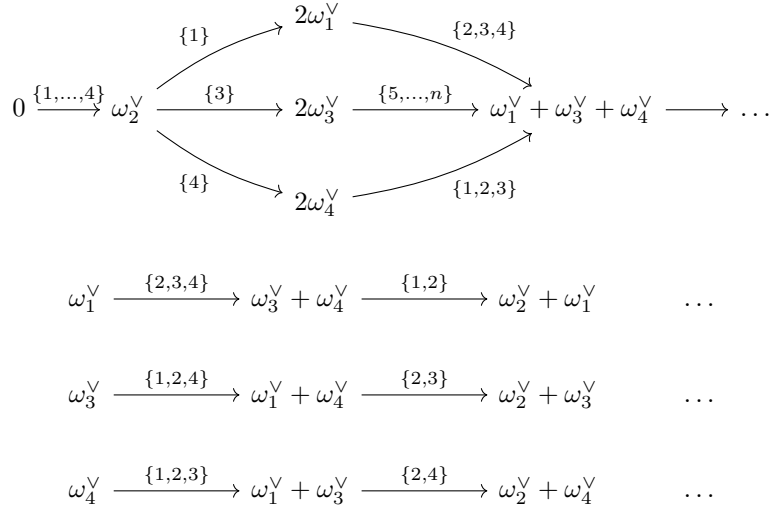
\begin{figure}
	\begin{tikzcd}
		 & && 2 \cw 1 \arrow[bend left=10]{drr}{\{2,3,4\}}&  \\
		0 \arrow{r}{\{1, \ldots, 4\}} & \cw 2 \arrow[bend left=10]{rru}{\{1\}} \arrow{rr}{\{3\}} \arrow[bend right=10]{drr}[swap]{\{4\}}&& 2 \cw 3  \arrow{rr}{\{5, \ldots, n\}} && \cw 1 + \cw 3 + \cw 4 \arrow{r}& \ldots \\
		&&& 2 \cw 4 \arrow[bend right=10]{rru}[swap]{\{1,2,3\}} 
	\end{tikzcd}
	\begin{tikzcd}
		&&&&\\
		\cw{1} \arrow{rr}{\{2,3,4 \}} && \cw 3 + \cw{4} \arrow{rr}{\{1,2\}}  && \cw 2 + \cw{1} & \ldots 
	\end{tikzcd}
	\begin{tikzcd}
		&&&&\\
		\cw{3} \arrow{rr}{\{1, 2,4 \}} && \cw 1 + \cw{4} \arrow{rr}{\{2, 3\}}  && \cw 2 + \cw{3} & \ldots 
	\end{tikzcd}
	\begin{tikzcd}
		&&&&\\
		\cw{4} \arrow{rr}{\{1, 2, 3\}} && \cw 1 + \cw{3} \arrow{rr}{\{2,4\}}  && \cw 2 + \cw{4} & \ldots 
	\end{tikzcd}
	
	\caption{The Bruhat partial order on dominant coweights of type $D_4$.}
	\label{figure-dom-cochars-D4}
\end{figure}

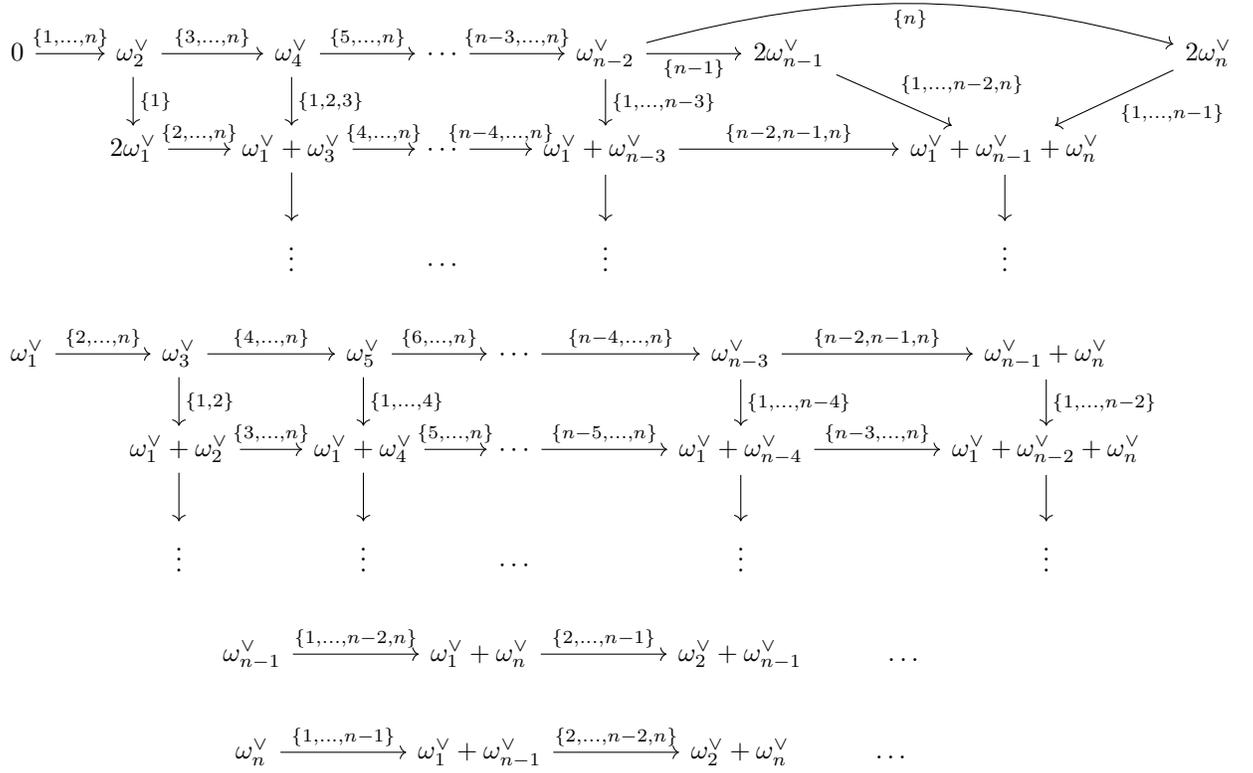
\begin{figure}
	\begin{tikzcd}
		0 \arrow{r}{\{1, \ldots, n\}} & \cw 2 \arrow{r}{\{3, \ldots, n\}} \arrow{d}{\{1\}}& \cw 4 \arrow{d}{\{1,2,3\}} \arrow{r}{\{5, \ldots, n\}} & \cdots \arrow{r}{\{n-3,\ldots,n\}}& \cw{n-2} \arrow{d}{\{1, \ldots,n-3\}}\arrow{r}[swap]{\{n-1\}} \arrow[bend left=15, swap]{rrr}{\{n\}} & 2\cw{n-1} \arrow{rd}{\{1,\ldots,n-2,n\}}&& 2\cw{n} \arrow{ld}{\{1, \ldots, n-1\}}   \\
		& 2 \cw 1 \arrow{r}{\{2, \ldots, n\}} & \cw 1+ \cw 3 \arrow{r}{\{4, \ldots, n\}} \arrow{d}  & \cdots \arrow{r}{\{n-4, \ldots, n\}} & \cw 1 + \cw{n-3} \arrow{rr}{\{n-2,n-1,n\}} \arrow{d} && \cw 1 + \cw{n-1} + \cw{n} \arrow{d} \\
		& & \vdots & \cdots & \vdots && \vdots & 
	\end{tikzcd}
	\begin{tikzcd}
		& \\
		\cw 1 \arrow{r}{\{2, \ldots, n\}} & \cw 3 \arrow{r}{\{4, \ldots, n\}} \arrow{d}{\{1,2\}}& \cw 5 \arrow{d}{\{1, \ldots, 4\}} \arrow{r}{\{6, \ldots, n\}} & \cdots \arrow{rr}{\{n-4, \ldots,n\}}&& \cw{n-3} \arrow{d}{\{1, \ldots,n-4\}}\arrow{rr}{\{n-2, n-1,n\}} && \cw{n-1} + \cw{n} \arrow{d}{\{1,\ldots,n-2\}}\\
		&  \cw 1 + \cw 2 \arrow{r}{\{3, \ldots, n\}} \arrow{d}& \cw 1+ \cw 4 \arrow{r}{\{5, \ldots, n\}} \arrow{d}  & \cdots \arrow{rr}{\{n-5, \ldots, n\}} && \cw 1 + \cw{n-4} \arrow{rr}{\{n-3,\ldots,n\}} \arrow{d} && \cw 1 + \cw{n-2} + \cw{n} \arrow{d} \\
		&\vdots & \vdots & \cdots && \vdots && \vdots  
	\end{tikzcd}

	\begin{tikzcd}
		&&&&\\
		\cw{n-1} \arrow{rr}{\{1, \ldots, n-2, n\}} && \cw 1 + \cw{n} \arrow{rr}{\{2, \ldots,n-1\}}  && \cw 2 + \cw{n-1} & \ldots 
	\end{tikzcd}
	\begin{tikzcd}
		&&&&\\
		\cw{n} \arrow{rr}{\{1, \ldots, n-1\}} && \cw 1 + \cw{n-1} \arrow{rr}{\{2, \ldots,n-2, n\}}  && \cw 2 + \cw{n} & \ldots 
	\end{tikzcd}
	
	\caption{The Bruhat partial order on dominant coweights of type $D_n$, $n\geq 6$ even.}
	\label{figure-dom-cochars-Dn-even}
\end{figure}

\begin{figure}
	\begin{tikzcd}
		0 \arrow{r}{\{1, \ldots, n\}} & \cw 2 \arrow{r}{\{3, \ldots, n\}} \arrow{d}{\{1\}}& \cw 4 \arrow{d}{\{1, \ldots, 3\}} \arrow{r}{\{5, \ldots, n\}} & \cdots \arrow{rr}{\{n-3, \ldots,n\}}&& \cw{n-3} \arrow{d}{\{1, \ldots,n-4\}}\arrow{rr}{\{n-2, n-1,n\}} && \cw{n-1} + \cw{n} \arrow{d}{\{1,\ldots,n-2\}}\\
		&  2\cw 1 \arrow{r}{\{2, \ldots, n\}} & \cw 1+ \cw 3 \arrow{r}{\{4, \ldots, n\}} \arrow{d}  & \cdots \arrow{rr}{\{n-5, \ldots, n\}} && \cw 1 + \cw{n-4} \arrow{rr}{\{n-3,\ldots,n\}} \arrow{d} && \cw 1 + \cw{n-2} + \cw{n} \arrow{d} \\
		& & \vdots & \cdots && \vdots && \vdots  
	\end{tikzcd}
	
	\begin{tikzcd}
		\cw 1 \arrow{r}{\{2, \ldots, n\}} & \cw 3 \arrow{r}{\{4, \ldots, n\}} \arrow{d}{\{1,2\}} 
		& \cdots \arrow{r}{\{n-3,\ldots,n\}}& \cw{n-2} \arrow{d}{\{1, \ldots,n-3\}}\arrow{r}[swap]{\{n-1\}} \arrow[bend left=15, swap]{rrr}{\{n\}} & 2\cw{n-1} \arrow{rd}{\{1,\ldots,n-2,n\}}&& 2\cw{n} \arrow{ld}{\{1, \ldots, n-1\}}   \\
		& \cw 1 + \cw 2 \arrow{r}{\{3, \ldots, n\}} \arrow{d} 
		& \cdots \arrow{r}{\{n-4, \ldots, n\}} & \cw 1 + \cw{n-3} \arrow{rr}{\{n-2,n-1,n\}} \arrow{d} && \cw 1 + \cw{n-1} + \cw{n} \arrow{d} \\
		&\vdots
		& \cdots & \vdots && \vdots & 
	\end{tikzcd}

	\begin{tikzcd}
		&&&&\\
		\cw{n-1} \arrow{rr}{\{1, \ldots, n-2, n\}} && \cw 1 + \cw{n} \arrow{rr}{\{2, \ldots,n-1\}}  && \cw 2 + \cw{n-1} & \ldots 
	\end{tikzcd}
	\begin{tikzcd}
		&&&&\\
		\cw{n} \arrow{rr}{\{1, \ldots, n-1\}} && \cw 1 + \cw{n-1} \arrow{rr}{\{2, \ldots,n-2, n\}}  && \cw 2 + \cw{n} & \ldots 
	\end{tikzcd}
	
	\caption{The Bruhat partial order on dominant coweights of type $D_n$, $n \geq 5$ odd.}
	\label{figure-dom-cochars-Dn-odd}
\end{figure}
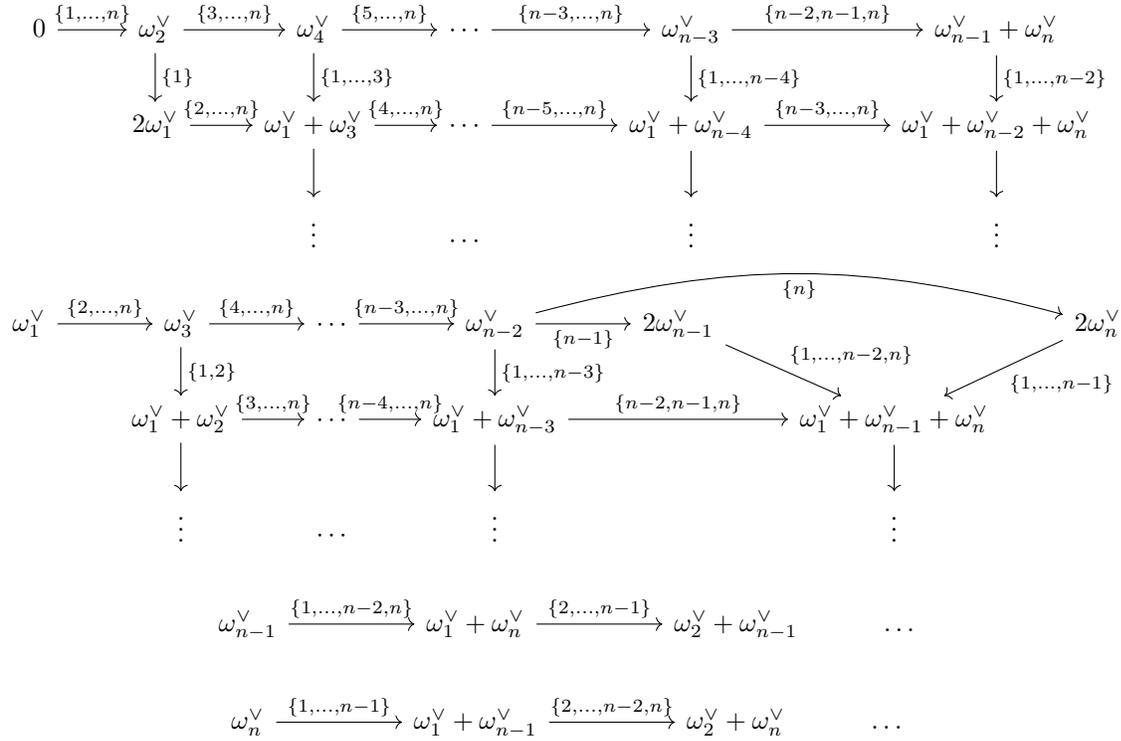

\bibliographystyle{alphaurl}
\bibliography{bib}

\end{document}